\documentclass[a4paper,11pt]{amsart}  

\usepackage{style}

\title{Lines on K3 quartic surfaces in characteristic 3}
\author{Davide Cesare Veniani}

\address{Institut für Algebraische Geometrie,
Leibniz Universität Hannover,
Wel\-fen\-garten 1, 30167 Hannover, Germany}
\curraddr{Institut für Geometrie und Topologie \\ Universität Stuttgart \\ Pfaffenwaldring 57 \\ 70569 Stuttgart, Germany}
\email{davide.veniani@mathematik.uni-stuttgart.de}

\keywords{quartic surface, K3 surface, line, positive characteristic, rational double point}

\thanks{The author acknowledges the financial support of the research training group GRK~1463 ``Analysis, Geometry and String Theory''.}

\date{December 14, 2020}

\subjclass[2010]{14J28, 14N10, 14N25}

\begin{document}

\begin{abstract}
We investigate the number of straight lines contained in a K3 quartic surface \(X\) defined over an algebraically closed field of characteristic~\(3\). We prove that if \(X\) contains \(112\) lines, then \(X\) is projectively equivalent to the Fermat quartic surface; otherwise, \(X\) contains at most \(67\) lines. We improve this bound to \(58\) if \(X\) contains a star (ie four distinct lines intersecting at a smooth point of \(X\)). Explicit equations of three \(1\)-dimensional families of smooth quartic surfaces with \(58\) lines, and of a quartic surface with \(8\) singular points and \(48\) lines are provided.
\end{abstract}

\maketitle

\section{Introduction}
A quartic surface in \(\IP^3\) whose singularities are isolated rational double points is called a \emph{K3 quartic surface}.
Together with \cite{veniani-char2} and \cite{veniani1}, this paper forms a trilogy dedicated to the study of the number of straight lines lying on a K3 quartic surface \(X\), denoted by \(|{\Fn X}|\).

The field \(\IK\) over which \(X\) is defined is always supposed to be algebraically closed.
The aim of this paper is to investigate the case \(\Char \IK = 3\). 

\subsection{Background}
Unlike smooth cubic surfaces, which always contain 27 lines, a dimension count shows that a general smooth quartic surface does not contain any line at all. One is therefore led to study the maximal number of lines that a quartic surface can contain.

Historically, the main focus has been on smooth complex quartic surfaces, which are an example of algebraic K3 surfaces. In 1882 F. Schur~\cite{schur} discovered a surface with 64 lines that now carries his name, given by the following equation:
\[
x^4 - xy^3 = w^4 - wz^3.
\]

The fact that \(|{\Fn X}| \leq 64\) for a smooth complex quartic surface \(X\) was proven by B. Segre in 1943~\cite{segre}. Around 70 years later, though, Rams and Schütt~\cite{rams.schuett:64.lines} discovered a flaw in Segre's argument (further investigated in a spin-off paper \cite{rams.schuett:quartics.lines.2nd.kind}) and fixed his proof, extending Segre's theorem to the case \(\Char \IK \neq 2,3\). 

Approximately at the same time, Degtyarev, Itenberg and Sertöz, spurred by  a remark of Barth~\cite{barth}, gave another solution with an algebraic approach. 
Their work resulted in a complete classification of large configurations of lines on smooth complex quartic surfaces up to projective equivalence~\cite{Degtyarev.Itenberg.Sertoz}. There are \(8\) possible configurations with more than \(52\) lines, corresponding to \(10\) distinct surfaces. A list of explicit equations for these surfaces, initiated by Schur \cite{schur}, Rams and Schütt~\cite{rams.schuett:64.lines}, Degtyarev, Itenberg and Sertöz~\cite{Degtyarev.Itenberg.Sertoz}, and Shimada and Shioda~\cite{shimada-shioda}, was completed by the author \cite{veniani:symmetries.equations}.

If \(\Char \IK \neq 2,3\), the inequality \(|{\Fn X}|\leq 64\) holds more generally for a (possibly non-smooth) K3 quartic surface \(X\) \cite{veniani1}. The bound of \(64\) lines is not known to be sharp for non-smooth K3 quartic surfaces. 
To the author's knowledge, the record for the number of lines contained in a non-smooth K3 quartic surface in the case \(\Char \IK = 0\) (respectively \(\Char \IK > 3\)) is attained by a surface with \(52\) lines (respectively by its reduction modulo \(5\) with \(56\) lines), whose equation is published in this paper for the first time (see \autoref{ex:52.lines.2.sing.pts}).

If \(\Char \IK = 2\), a sharp bound is known for both smooth and non-smooth K3 surfaces. More precisely, if \(X\) is smooth, then \(|{\Fn X}| \leq 60\). The bound was found by Degtyarev~\cite{degtyarev:lines.supersingular} and a surface attaining the bound was found by Rams and Schütt~\cite{rams.schuett:at.most.64.lines.char.2}. 
More generally, if \(X\) is a K3 quartic surface, \(|{\Fn X}| \leq 68\) and if equality holds, then \(X\) is not(!) smooth and projectively equivalent to a member of a \(1\)-dimensional family discovered by Rams and Schütt \cite{veniani-char2}.


\subsection{Principal results}
The Fermat quartic surface, which is defined by
\begin{equation} \label{eq:fermat} 
x^4+y^4+w^4+z^4 = 0,
\end{equation}
contains \(48\) lines over \(\IC\). Its reduction modulo \(3\) is a smooth, supersingular (ie of Picard rank~\(22\)) surface and contains \(112\) lines.
It is the smooth quartic surface with the highest number of lines in characteristic~\(3\), a result proven independently by Rams and Schütt \cite{RS-char3} and Degtyarev \cite{degtyarev:lines.supersingular}. 
The main result of this paper is the following generalization.

\begin{theorem}[see \autoref{subsec:proof.thm:char3}] \label{thm:char3}
Assume \(\Char \IK = 3\). If a K3 quartic surface \(X\) is not projectively equivalent to the Fermat surface, then \(|{\Fn X}| \leq 67\).
\end{theorem}

We refine \autoref{thm:char3} in the case that \(X\) contains four distinct lines meeting in a smooth point, a configuration which we call a \emph{star}. 
This configuration is particularly relevant in characteristic~\(3\). For instance, each line in the Fermat surface is contained in \(10\) stars, and if two lines on the Fermat surface intersect, then they belong to the same star.

\begin{addendum}[see \autoref{prop:star}] \label{add:char3}
Assume \(\Char \IK = 3\). If a K3 quartic surface \(X\) is not projectively equivalent to the Fermat surface and contains a star, then \(|{\Fn X}| \leq 58\).
\end{addendum}

For smooth surfaces more precise results were found by Degtyarev~\cite{degtyarev:lines.supersingular}, distinguishing between supersingular and non-supersingular surfaces.

A smooth supersingular quartic surface can contain \(112\), \(58\) or at most \(52\) lines (\cite[Theorem 1.2]{degtyarev:lines.supersingular}). 
There exist three different possible configurations of \(58\) lines. 
In \autoref{ex:58-1st-conf}, \autoref{ex:58-2nd-conf} and \autoref{ex:58-3rd-conf}, explicit equations of the three \(1\)-dimensional families of smooth quartic surfaces containing these configurations are provided. All of them contain a star. The first two families were discovered independently by Degtyarev (\cite[Proposition~8.13 and Proposition~8.14]{degtyarev:lines.supersingular}). To the author's knowledge, the third one is new. 


\autoref{add:char3} and the known results about smooth surfaces provide evidence for the following conjecture.

\begin{conjecture}
Assume \(\Char \IK = 3\). If a K3 quartic surface \(X\) is not projectively equivalent to the Fermat surface, then \(|{\Fn X}| \leq 58\).
\end{conjecture}

Neither  the bound of \(67\) lines in \autoref{thm:char3} nor the bound of \(58\) lines in \autoref{add:char3} are known to be sharp for non-smooth K3 quartic surfaces. In \autoref{ex:shimada-shioda} we show the existence of a non-smooth K3 quartic surface with \(48\) lines in characteristic~\(3\), the highest number known to the author. 

Finally, in \autoref{ex:52.lines.2.sing.pts} we exhibit a non-smooth K3 quartic surface with \(52\) lines in characteristic~\(0\) and \(56\) lines in characteristic~\(5\).

\subsection{Contents of the paper}
The main technique used in this paper is the study of the genus~\(1\) fibration induced by a line on a K3 quartic surface. A genus \(1\) fibration can be elliptic or a quasi-elliptic, the latter appearing only in characteristic 2 and 3 by a theorem of Tate~\cite{tate}. Heuristically,
quasi-elliptic fibrations can carry a higher number of fiber components than
elliptic fibrations, thus allowing for a higher number of lines on \(X\).
This technique is explained in \autoref{sec:K3-quartic}, where the notation is also fixed.

The technical core of the paper is contained in \autoref{sec:char3-valency}, where several bounds on the possible valencies of a line (see \autoref{def:valency}) are given.

The proofs of the main results are then carried out in~\autoref{sec:char3-proof}. 

Finally, various examples of K3 quartic surfaces with many lines are discussed in \autoref{sec:char3-examples}. 

\subsection*{Acknowledgments}
I am particularly indebted to Sławomir Rams and Matthias Schütt for suggesting the problem of studying the number of lines on K3 quartic surfaces. I also wish to thank Fabio Bernasconi, Alex Degtyarev and Víctor González-Alonso for many stimulating discussions.

\section{Preliminaries} \label{sec:K3-quartic}

In this section we fix our notation and recall some general facts about K3 quartic surfaces (see \cite[§2]{veniani-char2} for more details). 

\subsection{Lines and singular points}
Starting from here and throughout the paper, \(X\) denotes a K3 quartic surface defined over a field \(\IK\) of characteristic \(p \geq 0\), \(\Sing(X)\) is the set of singular points on \(X\) (which can be of type~\(\bA_n\), \(\bD_n\) or \(\bE_n\)), \(\rho : Z \rightarrow X\) is the minimal resolution of \(X\) (\(Z\) is a K3 surface), \(H\) is a hyperplane divisor in the linear system defined by \(\rho^*(\mathcal O_{X}(1))\), \(\ell\) is a line on \(X\), \(L\) is the strict transform of \(\ell\) on~\(Z\) and \(|{\Fn X}|\) is the number of lines lying on \(X\).

The pencil of planes \(\{\Pi_t\}_{t \in \IP^1}\) containing \(\ell\) induces a genus~1 fibration \(\pi: Z \rightarrow \IP^1\). A line \(\ell\) is \emph{elliptic} (respectively \emph{quasi-elliptic}) if it induces an elliptic (respectively \emph{quasi-elliptic}) fibration.
For \(t \in \IP^1\), the \emph{residual cubic} \(c_t\) is the union of the irreducible components of \(X \cap \Pi_t\) different from \(\ell\).
A fiber \(F_t\) of \(\pi\) is the pullback through \(\rho\) of the residual cubic \(c_t\) contained in \(X \cap \Pi_t\). A general fiber of \(\pi\) is denoted by \(F\).

The restriction of \(\pi\) to \(L\) is again denoted by \(\pi\). The \emph{degree} of \(\ell\) is the degree of the morphism \(\pi: L \rightarrow \IP^1\). (If \(\pi: L \rightarrow \IP^1\) is constant, we say that \(\ell\) has degree \(0\).)
The line \(\ell\) is \emph{separable} (respectively \emph{inseparable}) if \(\pi:L\rightarrow \IP^1\) is separable (respectively inseparable).
The \emph{singularity} of \(\ell\) is the cardinality of \(\Sing(X) \cap \ell\). 

A point \(P\) on a separable line \(\ell\) is a point \emph{of ramification \(n_l\)} if the corresponding point on \(L\) has ramification index \(n\) and \(\length(\Omega_{L/\IP^1}) = l\). If \(\Char \IK\) does not divide \(n\), then \(l = n-1\) and can be omitted, whereas if \(\Char \IK\) divides \(n\), then \(l \geq n\). We say that \(\ell\) has ramification \((n_l)^r(n'_{l'})^{r'}\ldots\) if it has \(r\) points of ramification \(n_l\) and so on.

\begin{lemma}[{\cite[Proposition 1.5]{veniani1}}] \label{lem:degree.singularity}
If \(\ell \subset X\) has degree~\(d\) and singularity \(s\), then \(d \leq 3 - s\), and \(d = 3\) if and only if \(s = 0\).
\end{lemma}

\begin{lemma}[{\cite[Lemma 2.7]{veniani1}}] \label{lemma:linesthroughsingularpoint}
    If \(P \in \Sing(X)\), then there are at most \(8\) lines on \(X\) passing through \(P\). \qed
\end{lemma}

\subsection{Valency}
Let \(\Pi \subset \IP^3\) be a plane such that \(X \cap \Pi\) splits into four lines \(\ell_1,\ldots,\ell_4\) (not necessarily distinct). If a line \(\ell' \subset X\) not lying on \(\Pi\) meets two or more distinct lines \(\ell_i\), then their point of intersection must be a singular point of the surface. It follows that \(|{\Fn X}|\) is bounded by
\begin{equation} \label{eq:FnX}
\begin{split}
    |{\Fn X}| &\leq \#\{\text{lines in \(\Pi\)}\} \\
            &\quad + \#\{\text{lines not in \(\Pi\) going through \(\Sing(X) \cap \Pi\)} \} \\ 
            &\quad + \sum_{i = 1}^4\#\{\text{lines not in \(\Pi\) meeting \(\ell_i\) in a smooth point}\}. \\
\end{split}
\end{equation}

\begin{definition} \label{def:valency}
The \emph{valency} of \(\ell\), denoted by \(v(\ell)\), is the number of lines \(\ell' \subset X\), \(\ell' \neq \ell\), such that \(\ell' \cap \ell \notin \Sing(X)\).
The \emph{local valency} of \(\ell\) at \(t \in \IP^1\), denoted by \(v_t(\ell)\), is the number of lines \(\ell'\subset c_t\), \(\ell'\neq \ell\), such that \(\ell' \cap \ell \notin \Sing(X)\).
\end{definition}

It holds
\begin{equation} \label{eq:v=sum-vF}
v(\ell) = \sum_{t \in \IP^1} v_t(\ell).
\end{equation}
A line \(\ell\) is said to be \emph{of type~\((p,q)\)} if there exist exactly \(p\) points \(t \in \IP^1\) such that \(c_t\) splits into three (not necessarily distinct) lines and exactly \(q\) points \(t \in \IP^1\) such that \(c_t\) splits into a line and an irreducible conic.
It follows from \eqref{eq:v=sum-vF} that if \(\ell\) has type~\((p,q)\) and degree~\(d\), then
\begin{equation} \label{eq:vl}
    v(\ell) \leq d p + q.
\end{equation}

\subsection{Lines of the first and second kind}

Let \(x_0,x_1,x_2,x_3\) be the coordinates of \(\IP^3\). We denote by \(\IK[x_0,\ldots,x_3]_{(d)}\) the space of forms (ie homogeneous polynomials) of degree~\(d\). Up to projective equivalence, we can suppose that the line \(\ell\) is given by \(x_0 = x_1 = 0\), so that the quartic \(X\) is defined by
\begin{equation} \label{eq:surfaceX}
 X: \sum_{i_0 + i_1 + i_2 + i_3 = 4} a_{i_0 i_1 i_2 i_3} x_{0}^{i_0} x_{1}^{i_1} x_{2}^{i_2} x_{3}^{i_3} = 0, \quad a_{i_0 i_1 0 0} = 0 \text{ for all \(i_0\), \(i_1\)},
\end{equation}
where \(i_0,\ldots,i_4\) are non-negative integers. 

The planes containing \(\ell\) are parametrized by \(\Pi_t: x_0 = t x_1\), where \(t = \infty\) denotes the plane \(x_1 = 0\). 
The residual cubic \(c_t\) is defined by the equation of \(\Pi_t\) and the equation \(g \in \IK[t][x_1,x_2,x_3]_{(3)}\) obtained by substituting \(x_0\) with \(t x_1\) in \eqref{eq:surfaceX} and factoring out \(x_1\). 
An explicit computation shows that there exist two forms \(\alpha,\beta \in \IK[x_2,x_3]_{(\deg(\ell))}\) such that \(\ell \cap c_t\) consists of the points \([0:0:x_2:x_3]\) satisfying
\begin{equation} \label{eq:penciloncurve}
 g(t;0,x_2,x_3) = t \alpha(x_2, x_3) + \beta(x_2, x_3) = 0.
\end{equation}

An \emph{inflection point} of a (possibly reducible) cubic curve \(c \subset \IP^2\) is a point of \(c\) which is a zero of the hessian of the cubic. If \(c\) contains a line as a component, all the points of the line are inflection points of \(c\). 
The following lemma can be checked by an explicit computation, too.

\begin{lemma} \label{lemma:line&tangentconic}
If \(c \subset \IP^2\) is a reducible cubic that is the union of an irreducible conic and a line \(m\), then the locus of inflection points of \(c\) is exactly~\(m\). \qed
\end{lemma}

If the surface \(X\) is defined as in \eqref{eq:surfaceX}, then the hessian of the equation \(g\) defining the residual cubic \(c_t\) restricted to the line \(\ell\) is given by 
\begin{equation} \label{eq:hessian}
h = \det \left( \frac{\partial^2 g}{\partial x_i x_j} \right)_{1\leq i,j, \leq 3} \bigg\rvert_{x_1 = 0} \in \IK[t][x_2,x_3]_{(3)},
\end{equation}
which is a polynomial of degree \(5\) in \(t\), with forms of degree \(3\) in \((x_2,x_3)\) as coefficients (if \(\Char \IK = 2\), this definition has to be slightly modified, see~\cite{RS-char3}).

The resultant \(R(\ell)\) with respect to the variable \(t\) of the polynomials \eqref{eq:penciloncurve} and \eqref{eq:hessian} is called the \emph{resultant} of the line \(\ell\).
We say that a line \(\ell\) of positive degree is a line of the \emph{second kind} if its resultant is identically equal to zero. Otherwise, we say that \(\ell\) is a line of the \emph{first kind}.

A root \([ x_2:  x_3]\) of \(R(\ell)\) corresponds to a point \(P = [0:0: x_2: x_3]\) on \(\ell\). If \(P\) is a smooth surface point, then it is an inflection point for the residual cubic passing through it.

\begin{proposition}[{\cite[Proposition 2.18]{veniani-char2}}] \label{prop:1stkind}
If \(\ell \subset X\) is a line of the first kind of degree~\(d\), then \(v(\ell) \leq 3 + 5\,d\). \qed
\end{proposition}

\subsection{Triangle-free surfaces} \label{subsec:triangle-free.K3}

We follow here the nomenclature of \cite[§5]{veniani1}. 
A Dynkin diagram (resp. extended Dynkin diagram) is also called an \emph{elliptic} graph (resp. \emph{parabolic} graph).
The \emph{line graph} \(\Gamma(X)\) of \(X\) is the dual graph of the strict transforms of its lines on \(Z\).
A \emph{triangle} on \(X\) is the union of three lines intersecting pairwise in smooth points of \(X\). We say that \(X\) is \emph{triangle-free} if it contains no triangles.

\begin{lemma}[{\cite[Lemma 2.22]{veniani-char2}}] \label{lem:elliptic-triangle-free}
If \(X\) is triangle-free and \(\ell \subset X\) is elliptic, then \(v(\ell) \leq 12\). \qed
\end{lemma}

\begin{lemma}[{\cite[Lemma 2.23]{veniani-char2}}] \label{lem:triangle-conf}
A triangle on \(X\) is contained in a plane~\(\Pi \subset \IP^3\) such that the intersection of \(\Pi\) and \(X\) is one of the configurations pictured in \autoref{fig:triangle-conf}. \qed
\end{lemma}

\begin{figure}[t]
\centering
\begin{tikzpicture}[line cap=round,line join=round,>=triangle 45,x=0.4cm,y=0.4cm,sing-pt/.style = {fill=white,thick}]
    \begin{scope}[xshift=0,yshift=0]
        \draw [name path=l1,thick] (-3,0) -- (3,0);
        \draw [name path=l2,thick] (-2.5,-0.87) -- (0.5,4.33);
        \draw [name path=l3,thick] (0,4.46) -- (0,-1);
        \draw [name path=l4,thick] (-2.13,2.23) -- (2.87,-0.5);
        
        \fill [name intersections={of=l1 and l2}] (intersection-1) circle (1.5pt);
        \fill [name intersections={of=l1 and l3}] (intersection-1) circle (1.5pt);
        \fill [name intersections={of=l1 and l4}] (intersection-1) circle (1.5pt);
        \fill [name intersections={of=l2 and l3}] (intersection-1) circle (1.5pt);
        \fill [name intersections={of=l2 and l4}] (intersection-1) circle (1.5pt);
        \fill [name intersections={of=l3 and l4}] (intersection-1) circle (1.5pt);

        \draw (2.5,3.5) node {\(\cA_0\)};
    \end{scope}
    \begin{scope}[xshift=80,yshift=0]
        \draw [name path=l1,thick] (-3,0) -- (3,0);
        \draw [name path=l2,thick] (-2.5,-0.87) -- (0.5,4.33);
        \draw [name path=l3,thick] (0,4.46) -- (0,-1);
        \draw [name path=l4,thick] (-2.13,2.23) -- (2.87,-0.5);
        
        \draw [name intersections={of=l1 and l2},sing-pt] (intersection-1) circle (2pt);
        \fill [name intersections={of=l1 and l3}] (intersection-1) circle (1.5pt);
        \fill [name intersections={of=l1 and l4}] (intersection-1) circle (1.5pt);
        \fill [name intersections={of=l2 and l3}] (intersection-1) circle (1.5pt);
        \fill [name intersections={of=l2 and l4}] (intersection-1) circle (1.5pt);
        \fill [name intersections={of=l3 and l4}] (intersection-1) circle (1.5pt);

        \draw (2.5,3.5) node {\(\cA_1\)};
    \end{scope}
    \begin{scope}[xshift=160,yshift=0]
        \draw [name path=l1,thick] (-3,0) -- (3,0);
        \draw [name path=l2,thick] (-2.5,-0.87) -- (0.5,4.33);
        \draw [name path=l3,thick] (0,4.46) -- (0,-1);
        \draw [name path=l4,thick] (-2.13,2.23) -- (2.87,-0.5);
        
        \draw [name intersections={of=l1 and l2},sing-pt] (intersection-1) circle (2pt);
        \draw [name intersections={of=l1 and l3},sing-pt] (intersection-1) circle (2pt);
        \fill [name intersections={of=l1 and l4}] (intersection-1) circle (1.5pt);
        \fill [name intersections={of=l2 and l3}] (intersection-1) circle (1.5pt);
        \fill [name intersections={of=l2 and l4}] (intersection-1) circle (1.5pt);
        \fill [name intersections={of=l3 and l4}] (intersection-1) circle (1.5pt);
        \draw (2.5,3.5) node {\(\cA_2\)};
    \end{scope}
    \begin{scope}[xshift=240,yshift=0]
        \draw [name path=l1,thick] (-3,0) -- (3,0);
        \draw [name path=l2,thick] (-2.5,-0.87) -- (0.5,4.33);
        \draw [name path=l3,thick] (0,4.46) -- (0,-1);
        \draw [name path=l4,thick] (-2.13,2.23) -- (2.87,-0.5);
        
        \draw [name intersections={of=l1 and l2},sing-pt] (intersection-1) circle (2pt);
        \draw [name intersections={of=l1 and l3},sing-pt] (intersection-1) circle (2pt);
        \draw [name intersections={of=l1 and l4},sing-pt] (intersection-1) circle (2pt);
        \fill [name intersections={of=l2 and l3}] (intersection-1) circle (1.5pt);
        \fill [name intersections={of=l2 and l4}] (intersection-1) circle (1.5pt);
        \fill [name intersections={of=l3 and l4}] (intersection-1) circle (1.5pt);
        \draw[color=black] (2.5,3.5) node {\(\cA_3\)};
    \end{scope}

    \begin{scope}[xshift=0,yshift=-80]
        \draw [name path=l1,thick] (-3,0)-- (3,0);
        \draw [name path=l2,thick] (-2.5,-0.87)-- (0.5,4.33);
        \draw [name path=l3,thick] (0,4.46)-- (0,-1);
        \draw [name path=l4,thick] (2.5,-0.87)-- (-0.5,4.33);

        \fill [name intersections={of=l1 and l2}] (intersection-1) circle (1.5pt);
        \fill [name intersections={of=l1 and l3}] (intersection-1) circle (1.5pt);
        \fill [name intersections={of=l1 and l4}] (intersection-1) circle (1.5pt);
        \fill [name intersections={of=l2 and l4}] (intersection-1) circle (1.5pt);
        
        \draw[color=black] (2.5,3.5) node {\(\cB_0\)};
    \end{scope}
    \begin{scope}[xshift=80,yshift=-80]
        \draw [name path=l1,thick] (-3,0)-- (3,0);
        \draw [name path=l2,thick] (-2.5,-0.87)-- (0.5,4.33);
        \draw [name path=l3,thick] (0,4.46)-- (0,-1);
        \draw [name path=l4,thick] (2.5,-0.87)-- (-0.5,4.33);

        \draw [name intersections={of=l1 and l2},sing-pt] (intersection-1) circle (2pt);
        \fill [name intersections={of=l1 and l3}] (intersection-1) circle (1.5pt);
        \fill [name intersections={of=l1 and l4}] (intersection-1) circle (1.5pt);
        \fill [name intersections={of=l2 and l4}] (intersection-1) circle (1.5pt);
        
        \draw (2.5,3.5) node {\(\cB_1\)};
    \end{scope}
    \begin{scope}[xshift=160,yshift=-80]
        \draw [name path=l1,thick] (-3,0)-- (3,0);
        \draw [name path=l2,thick] (-2.5,-0.87)-- (0.5,4.33);
        \draw [name path=l3,thick] (0,4.46)-- (0,-1);
        \draw [name path=l4,thick] (2.5,-0.87)-- (-0.5,4.33);

        \draw [name intersections={of=l1 and l2},sing-pt] (intersection-1) circle (2pt);
        \draw [name intersections={of=l1 and l3},sing-pt] (intersection-1) circle (2pt);
        \fill [name intersections={of=l1 and l4}] (intersection-1) circle (1.5pt);
        \fill [name intersections={of=l2 and l4}] (intersection-1) circle (1.5pt);
        
        \draw (2.5,3.5) node {\(\cB_2\)};
    \end{scope}
    \begin{scope}[xshift=240,yshift=-80]
        \draw [name path=l1,thick] (-3,0)-- (3,0);
        \draw [name path=l2,thick] (-2.5,-0.87)-- (0.5,4.33);
        \draw [name path=l3,thick] (0,4.46)-- (0,-1);
        \draw [name path=l4,thick] (2.5,-0.87)-- (-0.5,4.33);

        \draw [name intersections={of=l1 and l2},sing-pt] (intersection-1) circle (2pt);
        \draw [name intersections={of=l1 and l3},sing-pt] (intersection-1) circle (2pt);
        \draw [name intersections={of=l1 and l4},sing-pt] (intersection-1) circle (2pt);
        \fill [name intersections={of=l2 and l4}] (intersection-1) circle (1.5pt);
        \draw[color=black] (2.5,3.5) node {\(\cB_3\)};
    \end{scope}
    
    \begin{scope}[xshift=40,yshift=-160]
        \draw [thick] (-1.5, -0.6)-- (1.5, 4.6);
        \draw [thick] (-1.5,4.6)-- (1.5,-0.6);
        \draw [thick] (-3,2)-- (3,2);
        \draw [thick] (0,5)-- (0,-1);
        \fill (0,2) circle (1.5pt);
        \draw (2.5,3.5) node {\(\mathcal C_0\)};
    \end{scope}
    
    \begin{scope}[xshift=120,yshift=-160]
        \draw [name path=l2,thick] (-2.5,-0.87)-- (0.5,4.33);
        \draw [name path=l4,thick] (2.5,-0.87)-- (-0.5,4.33);
        \draw [thick] (-3,-0.1)-- (3,-0.1);
        \draw [thick] (-3,0.1)--(3,0.1);

        \fill (-2,0) circle (1.5pt);
        \fill (2,0) circle (1.5pt);
        \fill [name intersections={of=l2 and l4}] (intersection-1) circle (1.5pt);

        \draw [sing-pt] (-1,0) circle (2pt);
        \draw [sing-pt] (0,0) circle (2pt);
        \draw [sing-pt] (1,0) circle (2pt);

        \draw[color=black] (2.5,3.5) node {\(\mathcal D_0\)};
    \end{scope}
    \begin{scope}[xshift=200,yshift=-160]
        \draw [thick] (-2.5, -0.6)-- (0.5, 4.6);
        \draw [thick] (-2.5,4.6)-- (0.5,-0.6);
        \draw [thick] (-3,1.9)-- (3,1.9);
        \draw [thick] (-3,2.1)-- (3,2.1);
        
        \fill (-1,2) circle (1.5pt);

        \draw [sing-pt] (0,2) circle (2pt);
        \draw [sing-pt] (1,2) circle (2pt);
        \draw [sing-pt] (2,2) circle (2pt);
        
        \draw[color=black] (2.5,3.5) node {\(\mathcal E_0\)};
    \end{scope}
\end{tikzpicture}
\caption{Possible configurations of lines on a plane with a triangle. Singular points are marked white. In configurations \(\mathcal D_0\) and \(\mathcal E_0\) the singular points might coincide.}
\label{fig:triangle-conf}
\end{figure}

\subsection{Elliptic fibrations}
For further details on genus 1 fibration in any characteristic, we refer to \cite{bom-mum,cossec-dolgachev,rud-sha2,schuett.shioda:elliptic.surfaces}.

From now on, we suppose that \(\Char \IK =3\). We adopt Kodaira's notation for the possible types of a fiber \(F\) appearing in a genus \(1\) fibration. We denote by \(e(F)\) and \(\delta(F)\) respectively the Euler--Poincaré characteristic and the wild ramification index in characteristic~\(3\) of \(F\).
The possible values of \(e(F)\) and \(\delta(F)\) according to the type of \(F\) are displayed in \autoref{tab:fibers}.

The following formula holds for an elliptic fibration on a K3 surface.
\begin{equation} \label{eq:euler.elliptic}
    \sum_{t \in \IP^1} (e(F_t) + \delta(F_t)) = 24.   
\end{equation}

\begin{table}[b]
    \centering
    \caption{Euler--Poincaré characteristic and possible wild ramification indices of a fiber \(F\) in a genus 1 fibration.}
    \label{tab:fibers}
    \begin{tabular}{r|cccccccc}
    \toprule 
    type of \(F\) & \(\I_n\) & \(\II\) & \(\III\) & \(\IV\) & \(\I^*_n\) & \(\IV^*\) & \(\III^*\) & \(\II^*\) \\
    \midrule
    \(e(F)\) & \(n\) & \(2\) & \(3\) & \(4\) & \(n+2\) & \(8\) & \(9\) & \(10\) \\
    \(\delta(F)\) & \(0\) & \(\geq 1\) & \(0\) & \(\geq 1\) & \(0\) & \(\geq 1\) & \(0\) & \(\geq 1\) \\
    \bottomrule
    \end{tabular}
\end{table}

\subsection{Quasi-elliptic fibrations and cuspidal lines} 
A fiber \(F\) of a quasi-elliptic fibration can only have one of the following types: 
\(\II,\,\IV,\,\IV^*,\,\II^*\).

A K3 surface with a quasi-elliptic fibration is necessarily supersingular and the following formula holds
\begin{equation} \label{eq:euler.quasi-elliptic}
    \sum_{t \in \IP^1} (e(F_t) - 2) = 20.
\end{equation}

If \(iv\), \(iv^*\) and \(ii^*\) denote the numbers of reducible fibers of the respective type, then \eqref{eq:euler.quasi-elliptic} can be also written
\begin{equation} \label{eq:char3-euler}
iv + 3\,iv^* + 4\,ii^* = 10.
\end{equation}

\begin{figure}[t]
\centering
\begin{tikzpicture}[line cap=round,line join=round,>=triangle 45,x=0.75cm,y=0.75cm]
    \begin{scope}[xshift=0,yshift=0]
        \filldraw (-0.5,0) circle (2pt);
        \filldraw (0.5,0) circle (2pt);
        \filldraw (0,0.866) circle (2pt);
        \draw [thick] (-0.5,0)-- (0.5,0);

        \draw [thick] (0,0.866)--(-0.5,0);
        \draw [thick] (0,0.866)--(0.5,0);
        
        \draw [thick] (0,0.866)--(0,2);
        \draw [thick] (-0.5,0)--(0,2);
        \draw [thick] (0.5,0)--(0,2);
        \filldraw [fill=white,thick] (0,2) circle (3pt);
    \end{scope}

    \begin{scope}[xshift=100,yshift=0]
        \filldraw (-2,0) circle (2pt);
        \filldraw (-1,0) circle (2pt);
        \filldraw (0,0) circle (2pt);
        \filldraw (0,1) circle (2pt);
        \filldraw (0,2) circle (2pt);
        \filldraw (1,0) circle (2pt);
        \filldraw (2,0) circle (2pt);
        \draw [thick] (-2,0)-- (2,0);
        \draw [thick]  (0,0)-- (0,2);

        \draw [thick] (0,0)--(0,-1);
        \filldraw [fill=white,thick] (0,-1) circle (3pt);
    \end{scope}

    \begin{scope}[xshift=35,yshift=-70]
        \filldraw (-2,0) circle (2pt);
        \filldraw (-1,0) circle (2pt);
        \filldraw (0,0) circle (2pt);
        \filldraw (0,1) circle (2pt);
        \filldraw (1,0) circle (2pt);
        \filldraw (2,0) circle (2pt);
        \filldraw (3,0) circle (2pt);
        \filldraw (4,0) circle (2pt);
        \filldraw (5,0) circle (2pt);
        \draw [thick] (-2,0)-- (5,0);
        \draw [thick]  (0,0)-- (0,1);

        \draw [thick] (0,2)--(0,1);
        \filldraw [fill=white,thick] (0,2) circle (3pt);

        \draw [thick]  (3,1)-- (3,0);
        \filldraw [fill=white,thick]  (3,1) circle (3pt);

        \draw [thick]  (4.5,1)-- (4,0);
        \draw [thick]  (4.5,1)-- (5,0);
        \filldraw [fill=white,thick]  (4.5,1) circle (3pt);
    \end{scope}
\end{tikzpicture}
\caption{Cuspidal curve (white vertex) intersecting a reducible fiber of a quasi-elliptic fibration. Multiple white dots represent different possibilities.}
\label{fig:cuspidal.curve}
\end{figure}

The closure of the locus of cusps is a smooth curve~\(K\), called \emph{cuspidal curve}, such that \(K\cdot F = 3\). The restriction of the fibration to \(K\) is an inseparable morphism of degree \(3\). The cuspidal curve meets a reducible fiber in one of the ways picture in \autoref{fig:cuspidal.curve}.
In particular, the way \(K\) intersects \(F\) is uniquely determined unless \(F\) is of type~\(\II^*\). Note that \(K\) intersects a fiber of type~\(\IV\) at the intersection point of the three components.

A line \(\ell \subset X\) is said to be \emph{cuspidal} if it is quasi-elliptic and the cuspidal curve \(K\) of the induced fibration coincides with the strict transform \(L\) of~\(\ell\).
A cuspidal line is necessarily inseparable.

\section{Bounds on the valency} \label{sec:char3-valency}

From now on we assume that \(X\) is a K3 quartic surface defined over a field \(\IK\) of characteristic~3. This section is devoted to the study of the possible valencies of a line \(\ell \subset X\) (\autoref{def:valency}).

\subsection{Separable elliptic lines}
The results which we will prove presently are summarized in \autoref{tab:char3-elliptic}.

\begin{table}[b]
\caption{Bounds for the valency of a separable elliptic line.}
\centering
\label{tab:char3-elliptic}
\begin{tabular}{cccl} 
    \toprule
    kind & degree & singularity & valency \\
    \midrule
    \multirow{3}{*}{first kind}  & \(3\) & 0         & \(\leq 18\) (sharp bound) \\
                                  & \(2\) & 1         & \(\leq 13\) \\
                                  & \(1\) & 2 or 1    & \(\leq  8\) \\
    \midrule
    \multirow{4}{*}{second kind} & \(3\) & 0         & \(\leq 21\) (sharp bound) \\
                                  & \(2\) & 1         & \(\leq 14\) (sharp bound) \\
                                  & \(1\) & 2         & \(\leq 9\)  \\
                                  & \(1\) & 1         & \(\leq 11\) \\
    \midrule
    --                           & \(0\) & 3, 2 or 1 & \(\leq 2\) (sharp bound) \\
    \bottomrule
\end{tabular}
\end{table}

\begin{lemma} \label{lemma:2ndkind-ram2}
Let \(\ell \subset X\) be a separable line of the second kind and \(P \in \ell\) a smooth point of ramification~\(2\). Then, either the corresponding fiber is of type~\(\II\) with a cusp in \(P\), or the corresponding residual cubic splits into a double line and a simple line.
\end{lemma}
\proof Note that only lines of degree \(3\) and \(2\) can have a point \(P\) of ramification~\(2\).
We choose coordinates so that \(P\) is given by \([0:0:0:1]\). This means that
\[
a_{0103} = 0 \quad \text{and} \quad a_{0112} = 0.
\]
Since \(P\) is of ramification index~2 and it is non-singular, we can normalize
\[
a_{0121} = 1 \quad \text{and} \quad a_{1003} = 1.
\]
Since \(\ell\) is of the second kind, the following relations must be satisfied:
\[
a_{0202} = 0, \quad a_{0301} = a_{0211}^2 \quad \text{and} \quad a_{0310} = a_{0211}a_{0220}.
\]
This means that the residual cubic in \(x_0 = 0\) corresponding to \(P\) is given by
\[
\left(a_{0211} x_{1} - x_{2}\right)^2 x_{3} + f_3(x_1,x_2),
\]
where \(f_3\) is a form of degree 3. Looking at the explicit formula for \(f_3\) (which we do not write here), we see that either this cubic is irreducible and gives rise to a fiber of type~\(\II\), or the polynomial \(g = a_{0211} x_1 - x_2\) divides \(f_3\); in the latter case it is immediate to compute that also \(g^2\) divides \(f_3\).
\endproof

\begin{lemma} \label{lemma:2ndkind-ram3_4}
Let \(\ell \subset X\) be a separable line of the second kind and \(P \in \ell\) a point of ramification \(3_4\). Then, either the corresponding fiber is of type~\(\II\) with a cusp in \(P\), or the corresponding residual cubic splits into three lines (not necessarily distinct) intersecting in~\(P\).
\end{lemma}
\proof Note that \(\ell\) has necessarily degree~\(3\).
We choose coordinates so that \(P\) is given by \([0:0:0:1]\) and the fiber corresponds to the plane \(\Pi_0: x_0 = 0\). This means that
\[
a_{0103} = 0, \quad a_{0112} = 0 \quad \text{and} \quad a_{0121} = 0.
\]
A calculation with local parameters shows that \(\length(\Omega_{L/\IP^1}) = 4\) if and only if \(a_{1012} = 0\). 
Moreover, the following three coefficients must be different from 0: \(a_{0130}\), \(a_{1003}\) and \(a_{1021}\); 
the first two because otherwise there would be singular points on \(\ell\) (implying that the degree of \(\ell\) is less than \(3\)), the third because otherwise \(\ell\) would be inseparable. We can normalize them to 1, rescaling coordinates.
Two necessary conditions for the line \(\ell\) to be of the second kind are
\[
a_{0202} = 0 \quad \text{and} \quad a_{0211} = 0.
\]
Hence, the residual cubic in \(\Pi_0\) is given by
\[
a_{0301} x_{1}^{2} x_{3} + x_{2}^{3} + x_{1}f_2(x_1,x_2).
\]
Either the fiber is irreducible and has a cusp in \(P\) (\(a_{0301} \neq 0\)), or it splits into three concurrent lines (\(a_{0301} = 0\)).
\endproof

\begin{remark} \label{rmk:e=v}
Suppose that \(\ell\) is an elliptic line of degree \(3\). From~\autoref{tab:fibers} we observe 
\[
v_t(\ell) \leq e(F_t) \leq e(F_t) + \delta(F_t),
\]
since a reducible fiber has \(e(F_t) \geq 2\) and a fiber with at least three components has \(e(F_t) \geq 3\). From equations~\eqref{eq:v=sum-vF} and \eqref{eq:euler.elliptic} we infer that
\begin{equation} \label{eq:euler-valency}
    v(\ell) = \sum_{t\in\IP^1} v_t(\ell) \leq \sum_{t\in\IP^1} (e(F_t) + \delta(F_t)) = 24.
\end{equation}
If \(v_t(\ell) = e(F_t)\), then \(F_t\) is of type~\(\I_3\). Hence, if for any subset \(S \subset \IP^1\) one has
\[
\sum_{s\in S} e(F_s) = \sum_{s \in S} v_s(\ell) = N,
\]
for some \(N \geq 0\), then all fibers \(F_s\) must be of type~\(\I_3\) and, in particular, \(N\) must be divisible by~3.
\end{remark}

An application of the Riemann--Hurwitz formula yields the following lemma.

\begin{lemma}
If \(\ell \subset X\) is a separable line of degree~\(3\), then \(\ell\) has ramification \(2_1^4\), \(2_1 3_3\) or \(3_4\). \qed
\end{lemma}

\begin{proposition} \label{prop:lines2ndkind}
If \(\ell \subset X\) is a separable elliptic line of the second kind of degree~\(3\), then the valency of \(\ell\) is bounded according to \autoref{tab:lines2ndkind}.
\end{proposition}
\begin{table}[t]
\caption{Bounds for the valency of a separable elliptic line of the second kind of degree 3.}
\label{tab:lines2ndkind}
\centering
\begin{tabular}{cl} 
    \toprule
    ramification & valency \\
    \midrule
    \(2_1^4\)     & \(\leq 12\)     \\
    \(2_1 3_3\)    & \(\leq 21\) (sharp bound)  \\
    \(3_4\)       & \(\leq 21\) (sharp bound)   \\
    \bottomrule
\end{tabular}
\end{table}
\proof
Suppose first that \(\ell\) has a point \(P\) of ramification index \(2\).
According to \autoref{lemma:2ndkind-ram2}, the corresponding fiber \(F_P\) is either of type~\(\II\), so that \(e(F_P) + \delta_P \geq 2+1 = 3\) (\autoref{tab:fibers}) and \(v(F_P) = 0\), or it contains a double component, so that \(e(F_P) \geq 6\) and \(v(F_P) = 2\); in any case, the difference \(e(F_P) + \delta_P - v(F_P)\) is always at least 3. 
Therefore, if there are 4 points of ramification 2, then by formula~\eqref{eq:euler-valency} \(v(\ell)\) is not greater than \(24 - 4\cdot 3 = 12\), while if there is just one, \(v(\ell)\) is not greater than \(24 - 3 = 21\).

Suppose now that \(\ell\) has no point of ramification index \(2\), ie \(\ell\) has ramification~\(3_4\). 
If the ramified fiber \(F_0\) is of type~\(\II\), then there can be at most \(24-3 = 21\) lines meeting \(\ell\). Otherwise, \(F_0\) splits into three concurrent lines by \autoref{lemma:2ndkind-ram3_4}. Then, \(e(F_0) + \delta_0 \geq 5\) (type~\(\IV\) has wild ramification, too), which means that the contribution to \(v(\ell)\) of the other fibers is not greater than \(24-5 = 19\). Nonetheless, by \autoref{rmk:e=v} this contribution cannot be exactly \(19\), since \(19\) is not divisible by \(3\). Hence, again, we can have at most \(3 +18 = 21\) lines intersecting~\(\ell\).
\endproof

\begin{example}
The following surface contains a separable elliptic line (\(x_0 = x_1 = 0\)) of ramification \(2_1 3_3\) with valency \(21\): 
\[
    x_{0}^{4} + x_{0}^{2} x_{1} x_{2} - x_{1}^{3} x_{2} + x_{0} x_{1} x_{2}^{2} + x_{1} x_{2}^{3} + x_{0}^{2} x_{1} x_{3} + x_{1}^{2} x_{3}^{2} + x_{0} x_{2} x_{3}^{2} + x_{0} x_{3}^{3} = 0.
\]
\end{example}

\begin{example}
Let \(i\) be a square root of \(-1\). The following surface contains a separable elliptic line \(x_0 = x_1 = 0\) of ramification \(3_4\) with valency \(21\): 
\begin{multline*}
i x_{0}^{3} x_{1} + i x_{1}^{3} x_{2} + i x_{1} x_{2}^{3} - i x_{0}^{3} x_{3} + i x_{0} x_{1} x_{2} x_{3} + i x_{0} x_{3}^{3} \\
= x_{0}^{2} x_{1} x_{2} + x_{1}^{2} x_{2}^{2} + x_{0} x_{2}^{2} x_{3} - x_{0}^{2} x_{3}^{2}.
\end{multline*}
\end{example}

\begin{proposition} \label{prop:e-deg2}
If \(\ell \subset X\) is an elliptic line of degree~\(2\), then, \(v(\ell) \leq 14\).
\end{proposition}
\proof
By \autoref{prop:1stkind}, we can assume that \(\ell\) is of the second kind. Since \(\ell\) has degree \(2\), it must have singularity \(1\) (\autoref{lem:degree.singularity}). Let \(P\) be the singular point on~\(\ell\). 
The morphism \(\pi: L \rightarrow \IP^1\), being of degree 2, is separable and has two points of ramification index 2. At least one of the point of ramification must be different from \(P\). Let us call it \(Q\). By \autoref{lemma:2ndkind-ram2} either the fiber corresponding to \(Q\) is of type~\(\II\) or the residual cubic splits into a double line and a simple line.

Suppose the fiber \(F_Q\) is of type~\(\II\). If \(\ell\) is of type~\((p,q)\), then \(3\,p + 2\,q \leq 24-3=21\). Applying formula~\eqref{eq:vl}, we have \(v(\ell) \leq 2\,p + q = 14\).

If the residual cubic corresponding \(F_Q\) splits into a double line and a simple line, then it contributes 1 to the valency and at least 6 to the Euler--Poincaré characteristic. Applying formula~\eqref{eq:vl} again yields \(v(\ell) \leq 13\).
\endproof

\begin{example} 
The following surface contains an elliptic line \(\ell:x_0 = x_1 = 0\) of degree 2 with valency 14, thus attaining the bound in \autoref{prop:e-deg2}. The surface contains one point \(P = [0:0:0:1]\) of type~\(\bA_1\). The line \(\ell\) has 7 fibers of type~\(\I_3\) and one ramified fiber of type~\(\II\). The other ramified fiber corresponds to the plane \(x_0 = 0\) and is smooth. 
\[
x_{0}^{4} + x_{0}^{2} x_{1} x_{2} - x_{1}^{3} x_{2} + x_{0} x_{1} x_{2}^{2} + x_{1} x_{2}^{3} + x_{1}^{2} x_{3}^{2} + x_{0} x_{2} x_{3}^{2} = 0.
\]
\end{example}

\begin{proposition} \label{prop:e-deg1}
Let \(\ell \subset X\) be an elliptic line of degree \(1\). Then, \(v(\ell) \leq 9\) if \(\ell\) has singularity \(2\), and \(v(\ell) \leq 11\) if \(\ell\) has singularity \(1\).
\end{proposition}
\proof
The proof can be carried over word by word from the characteristic 0 case (see \cite[Propositions~2.13 and 2.14]{veniani1}).
\endproof

\subsection{Quasi-elliptic lines of degree 3}

\autoref{tab:char3-qe} summarizes the known bounds for the valency of a quasi-elliptic line, which we are about to prove.
\begin{table}[t]
\centering
\caption{Bounds for the valency of a quasi-elliptic line.}
\label{tab:char3-qe}
 \begin{tabular}{cll} \toprule
 degree &  & valency \\
 \midrule
 \multirow{2}{*}{\(3\)}  & cuspidal       & \(\leq 30\) (sharp bound) \\
                       & not cuspidal        & \(\leq 21\) (sharp bound) \\
          \(2\)          &        & \(\leq 14\) (sharp bound) \\
          \(1\)          &        & \(\leq 10\) \\
          \(0\)          &   & \(\leq 2\) \\
 \bottomrule
 \end{tabular}
\end{table}

\begin{lemma} \label{lemma:inseparable=>cuspidal}
If \(\ell \subset X\) is an inseparable line with \(v(\ell)>12\), then \(\ell\) is cuspidal (in particular, quasi-elliptic).
\end{lemma}
\proof
Up to coordinate change, we can suppose that the residual cubic contained in \(x_0 = t x_1\) intersects the line \(\ell:x_0 = x_1 = 0\) in \([0:0:0:1]\) for \(t = 0\) and in \([0:0:1:0]\) for \(t = \infty\). This means that the following coefficients vanish:
\[
a_{0103},\,a_{0112},\, a_{0121}; \, a_{1012},\,a_{1021},\, a_{1030}.
\]
Moreover, \(a_{1003}\) and \(a_{0130}\) must be different from \(0\), and can be normalized to \(1\) and \(-1\), respectively, by rescaling coordinates. Up to a Frobenius change of parameter \(t = s^3\), we can explicitly write the intersection point \(P_s\) of the residual cubic with \(\ell\), which is given by
\[
P_s = [0:0:s:1].
\]
If a residual cubic \(c_s\) is reducible, then all components must pass through \(P_s\); in particular, \(P_s\) must be a singular point of \(c_s\). One can see explicitly that \(P_s\) is a singular point of \(c_s\) if and only if \(s\) is a root of the following degree \(8\) polynomial:
\begin{equation} \label{eq:char3-phi}
\begin{split}
    g(s) & = a_{2020} s^{8} + a_{2011} s^{7} + a_{2002} s^{6} + a_{1120} s^{5} \\ 
& + a_{1111} s^{4} + a_{1102} s^{3} + a_{0220} s^{2} + a_{0211} s + a_{0202}.
\end{split}
\end{equation}
Furthermore, it can be checked by a local computation that if \(c_s\) splits into three (not necessarily distinct) lines, then \(s\) is a double root of \(g\). This implies that the valency of \(\ell\) is not greater than \(3\cdot 8 /2 = 12\), unless the polynomial \(g\) vanishes identically, but \(g \equiv 0\) implies that all points \(P_s\) are singular for \(c_s\), ie the line \(\ell\) is cuspidal.
\endproof

\begin{corollary} \label{cor:family-C}
If \(\ell\subset X\) is a cuspidal line, then \(X\) is projectively equivalent to a member of the family \(\mathcal C\) defined by
\[
\mathcal C\colon x_{0}x_{3}^3 - x_{1}x_{2}^3 + x_{2}q_3(x_{0},x_{1}) + x_{3}q_3'(x_{0},x_{1}) + q_4(x_{0},x_{1}),
\]
where \(q_3\), \(q_3'\) and \(q_4\) are forms of degree \(3\), \(3\) and \(4\), respectively.
\end{corollary}
\proof
The family can be found imposing that the polynomial \(g\) in the proof of \autoref{lemma:inseparable=>cuspidal} vanishes identically.
\endproof

\begin{corollary} \label{cor:cuspidal-reducible-fiber}
If \(\ell \subset X\) is a cuspidal line, then a residual cubic corresponding to a reducible fiber of \(\ell\) is either the union of three distinct concurrent lines, or a triple line.
\end{corollary}
\proof The intersection of a residual cubic with \(\ell\) is always one single point.

A residual cubic of \(\ell\) cannot be the union of a line and an irreducible conic, because the line and the conic would result in a fiber of type~\(\I_n\) (because the conic has to be tangent to \(\ell\)), which is impossible in a quasi-elliptic fibration. 

Therefore, a residual cubic relative to a reducible fiber must split into three (not necessarily distinct) lines. If at least two of them coincide, the plane on which they lie contains at least a singular point \(P\) of the surface (which is not on \(\ell\), since \(\ell\) has degree \(3\)). An explicit inspection of this configuration in the family \(\mathcal C\) (for instance, supposing up to change of coordinates that \(P\) is given by \([0:1:0:0]\)) shows that the residual cubic degenerates to a triple line.
\endproof

\begin{lemma} \label{lemma:qe-deg3}
If \(\ell \subset X\) is a quasi-elliptic line of degree \(3\), then \(v(\ell) \leq 30\).
\end{lemma}
\proof The fibration induced by the line \(\ell\) has at most 10 reducible fibers, each of which can contribute at most 3 to its valency.
\endproof

\begin{remark}
The bound of \autoref{lemma:qe-deg3} is sharp. As soon as a K3 quartic surface \(X\) is smooth, the valency of a quasi-elliptic line of degree \(3\) on \(X\) is automatically \(30\), because the fibration induced by \(\ell\) can only have \(10\) reducible fibers of type~\(\IV\), whose residual cubics are the union of three concurrent lines. Notably, this happens for all 112 lines of the Fermat surface.
\end{remark}

\begin{figure}[t]
\centering
\begin{tikzpicture}[line cap=round,line join=round,>=triangle 45,x=0.5cm,y=0.5cm,
sing-pt/.style={draw,thick,circle,fill=white,inner sep=1.5pt}]
    \begin{scope}[xshift=0,yshift=0]
        \draw [name path=res-cub,thick] (-2,1) .. controls (-1,-3) and (1.5,1) .. (1.5,3) ..controls (1.5,1) and (1.5,-1) .. (2.5,-0.5);
        \draw [name path=l,thick] (-3,0)-- (3,0);
            \node [above] at (3,0) {\(\ell\)};
        
        \node [sing-pt] at (1.5,3) [label=right: \(\bA_2\)] {};
        \fill [name intersections={of=res-cub and l}] (intersection-1) circle (2pt) (intersection-2) circle (2pt) (intersection-3) circle (2pt);
        
        \node at (0,-2) {\(\IV_{0}\)};
    \end{scope}
    
    \begin{scope}[xshift=100,yshift=0]
        \draw [name path=res-con,thick] (-2,2) arc (-180:70:1.5 and 2.5);
        \draw [name path=res-lin,thick] (1,4.33)-- (1,-1);
        \draw [name path=l,thick] (-3,0)-- (3,0);
            \node [above] at (3,0) {\(\ell\)};

        \node [sing-pt] at (1,2) [label=right: \(\bA_1\)] {};
        \fill [name intersections={of=res-con and l}] (intersection-1) circle (2pt) (intersection-2) circle (2pt);
        \fill [name intersections={of=res-lin and l}] (intersection-1) circle (2pt);
        
        \node at (0,-2) {\(\IV_{1}\)};
    \end{scope}
    
    \begin{scope}[xshift=200,yshift=0]
        \draw [name path=res-l2,thick] (-2.5,-0.87)-- (0.5,4.33);
        \draw [name path=res-l3,thick] (0,4.46)-- (0,-1);
        \draw [name path=res-l4,thick] (2.5,-0.87)-- (-0.5,4.33);
        \draw [name path=l,thick] (-3,0)-- (3,0);
            \node [above] at (3,0) {\(\ell\)};

        \fill [name intersections={of=res-l2 and l}] (intersection-1) circle (2pt);
        \fill [name intersections={of=res-l3 and l}] (intersection-1) circle (2pt);
        \fill [name intersections={of=res-l4 and l}] (intersection-1) circle (2pt);
        \fill [name intersections={of=res-l2 and res-l4}] (intersection-1) circle (2pt);

        \node at (0,-2) {\(\IV_{3}\)};
    \end{scope}
\end{tikzpicture}
\caption{Possible residual cubics corresponding to a fiber of type~\(\IV\).}
\label{fig:cubics-type-IV}
\end{figure}

\begin{lemma} \label{lemma:conf-IV}
Let \(\ell \subset X\) be any line of degree~\(3\). A fiber of type~\(\IV\) of the fibration induced by \(\ell\) corresponds to one of the residual cubics pictured in~\autoref{fig:cubics-type-IV}, with the only restriction that \(\ell\) cannot pass through a singular point.
\end{lemma}
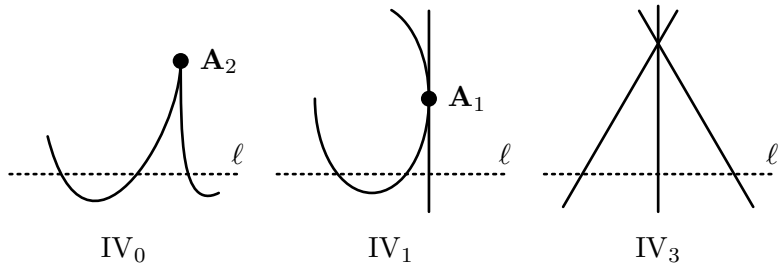
\proof
A fiber of type~\(\IV\) contains three simple components; hence, the corresponding residual cubic can also have only simple components. Since it cannot contain cycles, it must be one of the following, as in \autoref{fig:cubics-type-IV}: a cusp, a conic and a line meeting tangentially in one point, or three distinct lines meeting in one point.
The remaining components must come from the resolution of the singular points on the surface. The types of the singular points can be immediately deduced from the respective Dynkin diagrams.
\endproof

\begin{figure}[t]
\centering
\begin{tikzpicture}[line cap=round,line join=round,>=triangle 45,x=0.5cm,y=0.5cm,
sing-pt/.style={draw,thick,circle,fill=white,inner sep=1.5pt}]
    \begin{scope}[xshift=0,yshift=0]
        \draw [name path=res-cub,thick] (-2,1) .. controls (-1,-3) and (1.5,1) .. (1.5,3) .. controls (1.5,1) and (1.5,-1) .. (2.5,-0.5);
        \draw [name path=l,thick] (-3,0)-- (3,0);
            \node [above] at (3,0) {\(\ell\)};
        
        \node [sing-pt] at (1.5,3) [label=right: \(\bE_6\)] {};
        \fill [name intersections={of=res-cub and l}] (intersection-1) circle (2pt) (intersection-2) circle (2pt) (intersection-3) circle (2pt);
        
        \node at (0,-2) {\(\IV_{0}^*\)};
    \end{scope}
    
    \begin{scope}[xshift=100,yshift=0]
        \draw [name path=res-con,thick] (-2,2) arc (-180:70:1.5 and 2.5);
        \draw [name path=l,thick] (-3,0)-- (3,0);
            \node [above] at (3,0) {\(\ell\)};
        \draw [name path=res-lin,thick] (1,4.33)-- (1,-1);
        
        \node [sing-pt] at (1,2) [label=right: \(\bD_5\)] {};
        \fill [name intersections={of=res-con and l}] (intersection-1) circle (2pt) (intersection-2) circle (2pt);
        \fill [name intersections={of=res-lin and l}] (intersection-1) circle (2pt);

        \node at (0,-2) {\(\IV_{1}^*\)};
    \end{scope}
    
    \begin{scope}[xshift=200,yshift=0]
        \draw [thick] (-2.5,-0.87)-- (0.5,4.33);
        \draw [thick] (2.5,-0.87)-- (-0.5,4.33);
        \draw [thick] (-3,0)-- (3,0);
            \node [above] at (3,0) {\(\ell\)};
        \draw [thick] (0,4.46)-- (0,-1);
        \node [sing-pt] at (0,3.46) [label=right: \(\bD_4\)] {};
        \fill (-2,0) circle (2pt) 
              (0,0) circle (2pt) 
              (2,0) circle (2pt); 
        
        \node at (0,-2) {\(\IV_{3}^*\)};
    \end{scope}
    
    \begin{scope}[xshift=0,yshift=-110]
        \draw [name path=l2,thick] (2,4)-- (-2.5,-1);
        \draw [name path=l,thick] (-3,0)-- (3,0);
            \node [above] at (3,0) {\(\ell\)};
        \draw [thick] (1.1,4.33)-- (1.1,-1);
        \draw [thick] (0.9,4.33)-- (0.9,-1);
        \node [sing-pt] at (1,1.5) [label=right: \(\bA_5\)] {};
        \fill (1,0) circle (2pt) 
              (1,2.9) circle (2pt); 
        \fill [name intersections={of=l2 and l}] (intersection-1) circle (2pt);
        
        \node at (0,-2) {\(\IV_{2a}^*\)};
    \end{scope}
    \begin{scope}[xshift=100,yshift=-110]
        \draw [name path=l2,thick] (2,4)-- (-2.5,-1);
        \draw [name path=l,thick] (-3,0)-- (3,0);
            \node [above] at (3,0) {\(\ell\)};
        \draw [thick] (1.1,4.33)-- (1.1,-1);
        \draw [thick] (0.9,4.33)-- (0.9,-1);
        
        \node [sing-pt] at (1,2.9) [label=left: \(\bA_4\)] {};
        \node [sing-pt] at (1,1.5) [label=right: \(\bA_1\)] {};
        \fill [name intersections={of=l2 and l}] (intersection-1) circle (2pt) (1,0) circle (2pt);
        \node at (0,-2) {\(\IV_{2b}^*\)};
    \end{scope}
    \begin{scope}[xshift=200,yshift=-110]
        \draw [thick] (-3,0)-- (3,0);
            \node [above] at (3,0) {\(\ell\)};
        \draw [thick] (0.1,4.33)-- (0.1,-1);
        \draw [thick] (0.0,4.33)-- (0.0,-1);
        \draw [thick] (-0.1,4.33)-- (-0.1,-1);
        \node [sing-pt] at (0,3) [label=right: \(\bA_2\)] {};
        \node [sing-pt] at (0,2) [label=right: \(\bA_2\)] {};
        \node [sing-pt] at (0,1) [label=right: \(\bA_2\)] {};
        \fill (0,0) circle (2pt);
        
        \node at (0,-2) {\(\IV_{1t}^*\)};
    \end{scope}
\end{tikzpicture}
\caption{Possible residual cubics corresponding to a fiber of type~\(\IV^*\).}
\label{fig:cubics-type-IV*}
\end{figure}

\begin{lemma} \label{lemma:conf-IV*}
Let \(\ell \subset X\) be any line of degree~\(3\). A fiber of type~\(\IV^*\) of the fibration induced by \(\ell\) corresponds to one of the residual cubics pictured in \autoref{fig:cubics-type-IV*}, with the only restriction that \(\ell\) cannot pass through a singular point.
\end{lemma}
\proof
Besides the residual cubics with only simple components described in the previous lemma, we can also have multiple components, namely a double line and a simple line, or a triple line.

In the former case, the strict transforms of the lines can intersect (if their intersection point is smooth) or not (if their intersection point is singular), giving rise to two different configurations, which we distinguish by the letters \(a\) and \(b\). In the latter case, there is no ambiguity, since a fiber of type~\(\IV^*\) contains only one triple component.
\endproof

We denote by \(iv_0\), \(iv_1\),\ldots the number of fibers of type~\(\IV_0,\) \(\IV_1\), and so on. Note that the subscript indicates the local valency of the fiber.

\begin{lemma} \label{lem:char3-d=3-degK}
If \(\ell \subset X\) is a separable quasi-elliptic line of degree \(3\), then the degree of its cuspidal curve is at least \(3\) and at most \(7\). 
\end{lemma}
\proof Writing \(H = F + L\), one gets \(k = K\cdot H = 3 + K \cdot L\). The cuspidal curve \(K\) and the line \(L\) are distinct because \(\ell\) is separable.
The curve \(K\) can meet \(L\) only in points of ramification; moreover, a local computation shows that if \(K\) is tangent to \(L\), then ramification \(3_4\) occurs, and that higher order tangency cannot happen. We thus obtain the following bounds according to the ramification type of \(\ell\): \(K\cdot L \leq 4\) for type~\(2_1^4\), \(K\cdot L \leq 2\) for type~\(3_2 2_1\), and \(K\cdot L \leq 2\) for type~\(3_4\).
\endproof

\begin{proposition} \label{prop:qe-deg3}
If \(\ell \subset X\) is a separable quasi-elliptic line of degree \(3\), then \(v(\ell) \leq 21\).
\end{proposition}
\proof
A fiber of type~\(\II^*\) can have local valency at most \(2\), because it contains only one simple components and three distinct lines would give rise to three distinct simple components. Hence, recalling equation \eqref{eq:char3-euler},
\begin{align*}
v(\ell) & \leq 3\,iv + 3\,iv^* + 2\,ii^* \\
        & = 3\,(10-3\,iv^*-4\,ii^*) + 3\,iv^* + 2\,ii^* \\
        & = 30 - 6\,iv^* - 10\,ii^*.
\end{align*}
In particular, if \(ii^* > 0\), then \(v(\ell) \leq 20\), so we can suppose that \(\ell\) has no \(\II^*\)-fibers. Similarly, we can suppose that \(\ell\) has at most one fiber of type~\(\IV^*\).

If \(\ell\) has no \(\IV^*\)-fiber, then it must have \(10\) fibers of type~\(\IV\). Using the classification of \autoref{lemma:conf-IV}, we list the possible configurations with \(v(\ell)>21\) (16 cases) in \autoref{tab:char3-qe-deg3-v>16}.

\begin{table}[t]
\centering
\caption{Fibration types and valencies for a separable quasi-elliptic line \(\ell\) of degree 3 with \(v(\ell) \geq 16\).}
\label{tab:char3-qe-deg3-v>16}
\label{tab:qe-deg3}
 \begin{tabular}{cccccc} \toprule
 case & \(iv^*\) & \(iv_3\) & \(iv_1\) & \(iv_0\) & valency \\
 \midrule
  1 & -- & 10 & \(0\) & 0 & \(30\) \\
  2 & -- & 9 & \(1\) & 0 & \(28\) \\
  3 & -- & 9 & \(0\) & 1 & \(27\) \\
  4 & -- & 8 & \(2\) & 0 & \(26\) \\
  5 & -- & 8 & \(1\) & 1 & \(25\) \\
  6 & -- & 8 & \(0\) & 2 & \(24\) \\
  7 & -- & 7 & \(3\) & 0 & \(24\) \\
  8 & -- & 7 & \(2\) & 1 & \(23\) \\
  9 & -- & 7 & \(1\) & 2 & \(22\) \\
 10 & -- & 6 & \(4\) & 0 & \(22\) \\
 11 & \(iv^*_3\) &  7 & \(0\) & 0 & \(24\) \\
 12 & \(iv^*_3\) &  6 & \(1\) & 0 & \(22\) \\
 13 & \(iv^*_{2a}\) &  7 & \(0\) & 0 & \(23\) \\
 14 & \(iv^*_{2b}\) &  7 & \(0\) & 0 & \(23\) \\
 15 & \(iv^*_{1a}\) &  7 & \(0\) & 0 & \(22\) \\
 16 & \(iv^*_{1t}\) &  7 & \(0\) & 0 & \(22\) \\
 \bottomrule
 \end{tabular}
\end{table}

For each case, we consider the lattice generated by \(L\), a general fiber \(F\), the fiber components of the reducible fibers and the cuspidal curve \(K\) (which must be different from \(L\), since \(\ell\) is separable). All intersection numbers are univocally determined (\(L\cdot F = 3\) because \(\ell\) has degree~\(3\)), except for
\[
K \cdot L = K \cdot (H-F) = k - 3,
\]
but \(k\) can only take up the values \(3,\ldots,7\) on account of \autoref{lem:char3-d=3-degK}. We check that this lattice has rank bigger than \(22\) in all cases, except for case 6 with \(k=3\) (ie \(K \cdot L = 0\)).

On the other hand, this case does not exist. In fact, suppose that \(\ell\) is as in case 6 with \(K\cdot L = 0\); in particular, \(\ell\) has no ramified fibers with multiple components and, since \(v(\ell) = 24\), \(\ell\) is of the second kind. It follows that if \(\ell\) has a point of ramification \(2\), then by \autoref{lemma:2ndkind-ram2}, the ramified fiber must be a cusp, ie \(K\) intersects \(L\) so \(K\cdot L >0\);
if \(\ell\) has ramification \(3_4\), then by \autoref{lemma:2ndkind-ram3_4} the ramified fiber must be either a cusp or the union of three distinct lines; in both cases, \(K\cdot L >0\).
\endproof

\begin{example}
The following surface contains a separable quasi-elliptic line \(\ell: x_0 = x_1 = 0\) of degree~\(3\) with valency~21, thus attaining the bound of \autoref{prop:qe-deg3}:
\[
X: x_{1}^{4} + x_{0}^{2} x_{2}^{2} - x_{1}^{2} x_{2}^{2} - x_{1} x_{2}^{3} + x_{0} x_{2}^{2} x_{3} + x_{0} x_{3}^{3}.
\]
It contains only one singular point \([1:0:0:0]\) of type~\(\bE_6\).
\end{example}

\subsection{Quasi-elliptic lines of lower degree}

\begin{proposition} \label{prop:qe-deg2}
If \(\ell \subset X\) is a quasi-elliptic line of degree~\(2\), then \(v(\ell) \leq 14\).
\end{proposition}
\proof
By \autoref{prop:1stkind}, we can assume that \(\ell\) is of the second kind. Let \(P\) be the singular point on \(\ell\) and suppose \(t \in \IP^1\) is such that \(v_t(\ell) > 0\). The cubic \(c = c_t\) is reducible, because it contains at least a line. Suppose that \(c\) splits into a line \(m\) and an irreducible conic \(q\) (which must be tangent to each other because fibers of type~\(\I_n\) are not admitted in a quasi-elliptic fibration). Since \(v_t(\ell) >0\), the line \(m\) meets \(\ell\) in a smooth point; hence, the three configurations in \autoref{fig:quasi-elliptic.d=2.cubics} may arise.

\begin{figure}[b]
\centering
\begin{tikzpicture}[line cap=round,line join=round,>=triangle 45,x=0.5cm,y=0.5cm,
sing-pt/.style={draw,thick,circle,fill=white,inner sep=1.5pt}]
    \begin{scope}[xshift=0,yshift=0]
        \draw [name path=res-con,thick] (0,2.5) arc (-270:60:2 and 1.5);
        \draw [name path=l,thick] (-3,0)-- (3,0);
            \node [above] at (3,0) {\(\ell\)};
        \draw [name path=m,thick] (2,-1)-- (2,3.5);
        
        \node [sing-pt,label=below: \(P\)] at (-1.5,0) {};
        \fill (2,0) circle (2pt);
        \node [name intersections={of=m and res-con},sing-pt] at (intersection-1) {};
        \fill [name intersections={of=l and res-con}] (intersection-2) circle (2pt);
        \node [below] at (intersection-2) {\(Q\)};

        \node at (-2,2.5) {\(q\)};
        \node [right] at (2,3) {\(m\)};
    \end{scope}
    \begin{scope}[xshift=110,yshift=0]
        \draw [thick] (0,3) arc (-270:60:2 and 1.5);
        \draw [thick] (-3,0)-- (3,0);
            \node [above] at (3,0) {\(\ell\)};
        \draw [thick] (2,-1)-- (2,3.5);
        
        \node [sing-pt,label= below: \(P\)] at (0,0) {};
        \fill (2,0) circle (2pt);
        \fill [sing-pt] (2,1.5) circle (2pt);
        
        \node at (-2,3) {\(q\)};
        \node [right] at (2,3) {\(m\)};
    \end{scope}
    \begin{scope}[xshift=220,yshift=0]
        \draw [thick] (0,2.5) arc (-270:60:2 and 1.5);
        \draw [thick] (-3,1)-- (3,1);
            \node [above] at (3,1) {\(\ell\)};
        \draw [thick] (2,-1)-- (2,3.5);
        
        \node [sing-pt,label= below left: \(P\)] at (-2,1) {};
        \fill (2,1) circle (2pt);
        
        \node at (-2,2.5) {\(q\)};
        \node [right] at (2,3) {\(m\)};
    \end{scope}
\end{tikzpicture}
\caption{Possible residual cubics appearing in the proof of \autoref{prop:qe-deg2}.}
\label{fig:quasi-elliptic.d=2.cubics}
\end{figure}

The third configuration gives rise to a fiber of type~\(\III\), but this type does not appear in a quasi-elliptic fibration. (In the other configurations, the point of intersection of \(m\) with the residual conic is singular.) 
In the first configuration, if \(q \cap \ell = \{P,Q\}\), then \(Q\) is not an inflection point of the residual cubic by \autoref{lemma:line&tangentconic}, contradicting the fact that \(\ell\) is of the second kind. 
The second configuration can be ruled out by an explicit parametrization (in a line of the second kind, either the point \(P\) is a ramification point, or the cubic passing twice through \(P\) is singular at \(P\)).

Thus, all three configurations are impossible, so \(c\) must split into three (not necessarily distinct) lines and at least one of them should pass through \(P\). Since there can be at most 8 lines through a singular point (\autoref{lemma:linesthroughsingularpoint}), there can be at most 7 such reducible fibers, each of them contributing at most 2 to the valency of \(\ell\), whence \(v(\ell) \leq 14\).
\endproof

\begin{example} 
The bound given by \autoref{prop:qe-deg2} is sharp. In fact, the following quartic surface contains a quasi-elliptic line \(\ell: x_0 = x_1 = 0\) of degree \(2\) and valency \(14\):
\[
X : x_{0}^{4} + x_{0}^{3} x_{1} + x_{0} x_{1}^{3} + x_{1} x_{2}^{3} + x_{0} x_{1} x_{3}^{2} + x_{1}^{2} x_{3}^{2} + x_{0} x_{2} x_{3}^{2} = 0.
\]
The quartic contains two singular points, \(P = [0:0:0:1]\) of type~\(\bA_1\) and \(Q = [-1:1:1:0]\) of type~\(\bE_6\). The fibration induced by the line \(\ell\) has 7 fibers of type~\(\IV\) and one fiber of type~\(\IV^*\) corresponding to the plane containing~\(Q\). 
\end{example}

\begin{lemma}
If \(\ell \subset X\) is a quasi-elliptic line of degree \(1\), then \(v(\ell) \leq 10\).
\end{lemma}
\proof
The fibration induced by the line \(\ell\) has at most 10 reducible fibers, each of which contributes at most 1 to its valency.
\endproof

\section{Proof of \autoref{thm:char3}} \label{sec:char3-proof}

In order to prove the main theorem, we divide our arguments according to whether the K3 quartic surface contains a star (\autoref{subsec:star.case}) or more generally a triangle (\autoref{subsec:triangle.case}), or contains no triangles at all (\autoref{subsec:triangle-free.case}).
Throughout the section, \(X\) denotes a K3 quartic surface defined over a field of characteristic \(3\).

\subsection{Star case} \label{subsec:star.case}
We will start by analyzing quartic surfaces containing a star, which, we recall, is the union of four distinct lines meeting in a smooth point.

Since the four lines in a star are necessarily coplanar, a star is the same as a configuration~\(\mathcal C_0\) (see \autoref{fig:triangle-conf}). The lines have necessarily degree \(3\) because there are no singular points on the plane containing them.

We will first need a series of lemmas. In all of them, we will parametrize the surface \(X\) as in \eqref{eq:surfaceX} in such a way that the star is contained in the plane \(x_0 = 0\) and the lines meet at \([0:0:0:1]\), ie setting the following coefficients equal to \(0\):
\[
a_{0301},\,a_{0211},\,a_{0121},\,a_{0202},\,a_{0112},\,a_{0103}.
\]
If necessary, we will parametrize a second line in the star \(\ell'\) as \(x_0 = x_2 = 0\), by further assuming \(a_{0400} = 0\).

\begin{lemma} \label{lemma:star-1st-kind}
If \(\ell \subset X\) is a line of the first kind in a star, then \(v(\ell) \leq 15\).
\end{lemma}
\proof
It can be checked by an explicit computation that the resultant of \(\ell\) has a root of order \(6\) at the center of the star; this implies that there are at most \(18-6 = 12\) lines meeting \(\ell\) not contained in the star.
\endproof

\begin{lemma} \label{lemma:star-2_1 3_3}
If \(\ell \subset X\) is a separable line of ramification \(2_1 3_3\) contained in a star, then it is of the first kind.
\end{lemma}
\proof
By a change of coordinates, we can assume that the point of ramification index 2 is \([0:0:1:0]\), and that ramification occurs at \(x_1 = 0\). Imposing that \(\ell\) is of the second kind leads to a contradiction (\(\ell\) cannot be separable).
\endproof

\begin{lemma} \label{lemma:star-3_4}
If three lines in a star contained in \(X\) are separable and at least two of them have ramification \(3_4\), then the third one also has ramification \(3_4\).
\end{lemma}
\proof Beside \(\ell:x_0 = x_1 = 0\) and \(\ell':x_0 = x_2 = 0\), we can suppose without loss of generality that a third line is given by \(\ell'':x_0 = x_1 + x_2 = 0\), setting \(a_{0220} = a_{0130}+a_{0310}\). The conditions for \(\ell\), \(\ell'\) or \(\ell''\) to be of ramification~\(3_4\) are \(a_{1012} = 0\), \(a_{1102} = 0\) and \(a_{1012} = a_{1102}\), respectively. Two of them imply the third one.
\endproof

\begin{lemma} \label{lemma:star-no-3x3_4}
If three lines in a star contained in \(X\) are separable, then at most two of them can be of the second kind.
\end{lemma}
\proof
We parametrize \(\ell\), \(\ell'\) and \(\ell''\) as in the previous Lemma. Imposing that all three of them are of the second kind leads to a contradiction (at least one of them must be inseparable).
\endproof

\begin{lemma} \label{lemma:star-1st-kind+3_4}
Let \(\ell, \ell' \subset X\) be two lines in a star. If \(\ell\) is a separable line of the second kind, and \(\ell'\) is a line of the first kind of ramification~\(3_4\), then \(v(\ell')\leq 12\).
\end{lemma}
\proof
This can be checked again by an explicit computation of the resultant of \(\ell'\), which has now a root of order \(9\) at the center of the star.
\endproof

\begin{lemma} \label{lemma:star-not-cuspidal+cuspidal}
Let \(\ell, \ell' \subset X\) be two lines in a star. If \(\ell\) is a cuspidal line,
and \(\ell'\) is not cuspidal, then \(v(\ell')\leq 12\).
\end{lemma}
\proof
We parametrize \(\ell: x_0 = x_1 = 0\) as in \autoref{cor:family-C}. By virtue of \autoref{lemma:inseparable=>cuspidal}, we can suppose that \(\ell':x_0 = x_2 = 0\) is separable. An explicit computation shows that \(\ell'\) cannot be of the second kind, and that its resultant has a root of order \(9\) in \(x_2 = 0\).
\endproof

\begin{lemma} \label{lemma:star-not-cuspidal+2xcuspidal}
Let \(\ell, \ell', \ell'' \subset X\) be three lines in a star. If \(\ell\) and \(\ell'\) are cuspidal, and \(\ell''\) is not cuspidal, then \(v(\ell'') = 3\).
\end{lemma}
\proof
We parametrize \(\ell:x_0 = x_1 = 0\) and \(\ell':x_0 = x_2 = 0\) as in \autoref{cor:family-C}, ie we suppose that \(X\) is given by the family \(\mathcal C\) where the following coefficients are set to zero:
\[
a_{0400},\,a_{0301};\,a_{1201},\,a_{1300};\,a_{2200},\,a_{2101},\,a_{1210}.
\]
By a further rescaling we put \(a_{0310} = 1\) and we consider \(\ell'': x_0 = x_1-x_2 = 0\). The line~\(\ell''\) is inseparable and we can compute its polynomial \(g\) as in formula \eqref{eq:char3-phi} in the proof of \autoref{lemma:inseparable=>cuspidal} (by parametrizing the pencil with \(x_0 = s^3(x_1-x_2)\)), which turns out to be
\[
g(s) = a_{2110} s^8.
\]
This means that \(\ell''\) has only one singular fiber in \(s = 0\) (namely a fiber of type~\(\IV\) with the maximum possible index of wild ramification), unless \(a_{2110} = 0\) and \(g \equiv 0\), in which case \(\ell''\) is cuspidal.
\endproof

\begin{lemma} \label{lemma:cuspidal=>two-stars}
If \(\ell \subset X\) is a cuspidal line which is not contained in at least two stars, then \(v(\ell) \leq 6\).
\end{lemma}
\proof
On account of \autoref{cor:cuspidal-reducible-fiber}, the number of stars in which \(\ell\) is contained is exactly equal to the number of fibers of type~\(\IV\) in its fibration; moreover, \(v_t(\ell) = 1\) if \(t \in \IP^1\) corresponds to a fiber of type~\(\IV^*\) or \(\II^*\), yielding 
\[ 
    v(\ell) = 3\,iv + iv^* + ii^*.
\]
Recalling formula \eqref{eq:char3-euler}, we deduce that if \(iv < 2\) then \(v(\ell) \leq 6\).
\endproof

\begin{proposition} \label{prop:star}
If \(X\) contains a star and is not projectively equivalent to the Fermat surface, then \(|{\Fn X}| \leq 58\).
\end{proposition}
\proof
If \(\ell_1,\,\ell_2,\,\ell_3,\,\ell_4\) are the lines contained in the star, it holds from \eqref{eq:FnX}
\[
|{\Fn X}| \leq 4 + \sum_{i = 1}^4 (v(\ell_i)-3) = \sum_{i=1}^4 v(\ell_i) -8.
\]

Suppose first that all lines \(\ell_i\) are not cuspidal.
If \(v(\ell_i) \leq 15\) for \(i=1,2,3,4\), then 
\[ 
    |{\Fn X}| \leq 4\cdot 15-8 = 52.
\]

If \(v(\ell_1)>15\), then by \autoref{lemma:star-1st-kind} and \autoref{lemma:inseparable=>cuspidal}, \(\ell_1\) must be separable of the second kind; hence \(v(\ell_1) \leq 21\); if \(v(l_i)\leq 15\) for \(i=2,3,4\), then 
\[
    |{\Fn X}| \leq (21+3\cdot 15)-8 = 58.
\]
If \(v(\ell_1)>15\) and \(v(\ell_2)>15\), then both \(\ell_1\) and \(\ell_2\) are separable lines of the second kind. On account of \autoref{lemma:star-2_1 3_3}, they both have ramification \(3_4\). We claim that both \(v(\ell_3)\) and \(v(\ell_4)\) are not greater than \(12\). In fact, if \(\ell_3\) is separable, then by \autoref{lemma:star-3_4} and \autoref{lemma:star-no-3x3_4} it must be of the first kind and have ramification \(3_4\), which in turn implies that \(v(\ell_3)\leq 12\), because of \autoref{lemma:star-1st-kind+3_4}; if \(\ell_3\) is inseparable, then \(v(\ell) \leq 12\) by \autoref{lemma:inseparable=>cuspidal}. The same applies to \(\ell_4\). Hence, we conclude that \[
    |{\Fn X}|\leq (2\cdot 21 + 2\cdot 12) -8 = 58.
\]

Assume now that exactly one of the lines, say \(\ell_1\), is cuspidal, so that \(v(\ell_1) \leq 30\). On account of \autoref{lemma:star-not-cuspidal+cuspidal} we have \[|{\Fn X}| \leq (30 + 3\cdot 12) -8 = 58.\]

Suppose then that both \(\ell_1\) and \(\ell_2\) are cuspidal. If \(\ell_3\) and \(\ell_4\) are not cuspidal, then by \autoref{lemma:star-not-cuspidal+2xcuspidal} \[|{\Fn X}| \leq (2\cdot 30 + 2\cdot 3) - 8 = 58.\]

Finally, suppose that \(\ell_1\), \(\ell_2\) and \(\ell_3\) are cuspidal.
By a local computation it can be seen that \(\ell_4\) is also necessarily cuspidal.
Thanks to the bound of \autoref{lemma:cuspidal=>two-stars}, we can suppose that at least two lines, say \(\ell_1\) and \(\ell_2\), are part of another star.
Pick two lines \(\ell'_1\) and \(\ell'_2\), each of them in another star containing \(\ell_1\) respectively \(\ell_2\), which intersect each other (necessarily in a smooth point).
Perform a change of coordinates so that \(\ell_1\), \(\ell_2\), \(\ell'_1\) and \(\ell'_2\) are given respectively by \(x_0 = x_1 = 0\), \(x_0 = x_2 = 0\), \(x_1 = x_3 = 0\) and \(x_2 = x_3 = 0\).
Impose that \(\ell_1\), \(\ell_2\) and \(\ell_3\) are cuspidal lines: the resulting surface is projectively equivalent to Fermat surface.
\endproof

\subsection{Triangle case} \label{subsec:triangle.case}

In this section we study the case in which \(X\) admits a triangle. The three lines forming the triangle need to be coplanar, and we denote by \(\Pi\) the plane on which they lie. 
The plane \(\Pi\) intersects \(X\) also in a fourth line, which might coincide with one of the first three. 

\begin{proposition} \label{prop:triangle-mult}
If \(X\) admits a configuration \(\mathcal D_0\) or \(\mathcal E_0\), then \(|{\Fn X}| \leq 60\).
\end{proposition}
\proof 
Let \(\ell_0\) be the double line in the plane \(\Pi\) containing one of the two configurations, and let \(\ell_1\) and \(\ell_2\) be the two simple lines. Lines meeting \(\ell_0\) different from \(\ell_1\) and \(\ell_2\) must pass through the singular points; hence, by \autoref{lemma:linesthroughsingularpoint} there can be at most \(3\cdot(8-1) = 21\) of them.

Note that \(\ell_1\) and \(\ell_2\) cannot be cuspidal because of \autoref{cor:cuspidal-reducible-fiber}.

In the fibrations induced by \(\ell_1\) and \(\ell_2\) the plane \(\Pi\) corresponds to a fiber with a multiple component, hence with Euler--Poincaré characteristic at least 6; therefore, if \(\ell_1\) and \(\ell_2\) are both elliptic, there can be at most \(18\) more lines meeting them, yielding by \eqref{eq:FnX}
\[
|{\Fn X}| \leq 3\cdot(8-1) + (18 + 18) + 3 = 60.
\]

Suppose that \(\ell_1\) is quasi-elliptic. The plane \(\Pi\) corresponds to a fiber of type~\(\IV^*\) or~\(\II^*\); hence, there can be only one singular point on \(\ell_0\): in fact, by inspection of the Dynkin diagrams, a component of multiplicity \(2\) in these fiber types meets at most \(2\) other components (and one of them is the strict transform of \(\ell_2\)). The lines \(\ell_1\) and \(\ell_2\) not being cuspidal, we know that they have valency at most \(21\). It follows from \eqref{eq:FnX} that
\[
 |{\Fn X}| \leq    (8-1) + 2\cdot(21-2) + 3 = 48. \qedhere
\]
\endproof

\begin{lemma} \label{lemma:A0A1}
Let \(\ell, \ell' \subset X\) be two lines of degree \(3\) in configuration \(\cA_0\) or \(\cA_1\). If \(v(\ell)>18\), then \(v(\ell')\leq 18\).
\end{lemma}
\proof
Let \(\Pi\) be the plane containing \(\ell\) and \(\ell'\). Both lines are separable, since otherwise the respective residual cubics would intersect them in one point. We suppose that also \(v(\ell')>18\) and look for a contradiction.

Since both lines have valency greater than \(18\), they must be lines of the second kind with ramification \(2_1 3_3\) or \(3_4\).
In particular, they must have a point of ramification \(3\) (let us call it \(P\in\ell\) and \(P'\in\ell'\)), which does not lie on \(\Pi\).
Up to change of coordinates, we can assume the following:
\begin{enumerate}
\item \(\Pi\) is the plane \(x_0 = 0\);
\item \(\ell\) and \(\ell'\) are given respectively by \(x_0 = x_1 = 0\) and \(x_0 = x_2 = 0\);
\item \(P\) is given by \([0:0:1:0]\) and \(P'\) by \([0:1:0:0]\);
\item ramification in \(P\) (resp. \(P'\)) occurs in \(x_1 = 0\) (resp. \(x_2 = 0\)).
\end{enumerate}
This amounts to setting the following coefficients equal to \(0\):
\[
a_{0400},\,a_{0301},\, a_{0202},\, a_{0103};\, a_{1030},\, a_{1021},\, a_{1012};\, a_{1300},\, a_{1201},\, a_{1102}.
\]
Furthermore, \(a_{0112} \neq 0\), since the two residual lines in \(\Pi\) do not contain \([0:0:0:1]\), the intersection point of \(\ell\) and \(\ell'\); we set \(a_{0112} = 1\) after rescaling one variable.

Two necessary condition for \(\ell\) and \(\ell'\) to be lines of the second kind are
\[
a_{0310} = a_{0211}^2
\quad \text{and} \quad a_{0130} = a_{0121}^2. \quad
\]
This means that the residual conic in \(\Pi: x_0 = 0\) is given explicitly by
\begin{equation} \label{eq:Q0}
a_{0211}^{2} x_{1}^{2} + a_{0121}^{2} x_{2}^{2} + a_{0220} x_{1} x_{2} + a_{0211} x_{1} x_{3} + a_{0121} x_{2} x_{3} + x_{3}^{2} = 0.
\end{equation}
This conic splits into two lines by hypothesis; hence, it has a singular point. Computing the derivatives, one finds that the following condition must be satisfied:
\[
a_{0220} = -a_{0121} a_{0211}.
\]
Substituting into \eqref{eq:Q0}, one finds that the conic degenerates to a double line:
\[
(a_{0211} x_{1} + a_{0121} x_{2} - x_{3})^2 = 0;
\]
thus, we have neither configuration \(\cA_0\) nor \(\cA_1\).
\endproof

\begin{proposition} \label{prop:triangle-not-star}
If \(X\) contains a triangle but not a star, then \(|{\Fn X}| \leq 67\).
\end{proposition}
\proof
The proof is a case-by-case analysis on the configurations that are given by \autoref{lem:triangle-conf}, except configurations \(\mathcal C\) (a star, treated in \autoref{prop:star}), \(\mathcal D_0\) and \(\mathcal E_0\) (treated in \autoref{prop:triangle-mult}). 

We recall the bounds on the valency of \autoref{sec:char3-valency} and the fact that there are at most \(8\) lines through a singular point.

Note that all lines in a configuration of type~\(\cA\) are elliptic, since the configuration corresponds to a fiber of type~\(\I_n\).

Suppose that \(X\) contains configuration \(\cA_0\). If one of the four lines has valency \(v(\ell) > 18\), then the valency of all other three lines is not greater than \(18\) by \autoref{lemma:A0A1}. It follows from \eqref{eq:FnX} that
\[
    |{\Fn X}| \leq 4 + (21-3) + 3\cdot (18-3) = 67.
\]

Suppose that \(X\) contains configuration \(\cA_1\). Let \(\ell_1\) and \(\ell_2\) be the lines through the singular point, and \(\ell_3\) and \(\ell_4\) the other two lines.
We know that \(v(\ell_i) \leq 14\), \(i=1,2\), whereas \autoref{lemma:A0A1} applies to \(\ell_3\) and \(\ell_4\), yielding \(v(\ell_3) + v(\ell_4) \leq 18+21\).
It follows from \eqref{eq:FnX} that
\[
    |{\Fn X}| \leq 4 + (8-2) + 2\cdot (14-2) + (18-3) + (21-3) = 67.
\]

Suppose that \(X\) contains configuration \(\cA_2\). Let \(\ell_1\) be the line of singularity~\(2\). According to \autoref{prop:e-deg1}, \(v(\ell_1) \leq 9\). Therefore, it holds
\[
    |{\Fn X}| \leq 4 + 2\cdot (8-2) + (9-1) + 2\cdot (14-2) + (21-3) = 66.
\]

If \(X\) contains configuration \(\cA_3\), then 
\[
    |{\Fn X}| \leq 4 + 3\cdot (8-2) + 2 + 3\cdot (14-2) = 60.
\]

If \(X\) contains a configuration of type~\(\cB\), then the three lines meeting at the same (smooth) point must be of the first kind by \autoref{lemma:2ndkind-ram2}. It follows from \eqref{eq:FnX} that
\[
    |{\Fn X}| \leq \begin{cases}
    67 & \text{if \(X\) contains configuration \(\cB_0\)}, \\
    63 & \text{if \(X\) contains configuration \(\cB_1\)}, \\
    61 & \text{if \(X\) contains configuration \(\cB_2\)}, \\
    57 & \text{if \(X\) contains configuration \(\cB_3\)}.
    \end{cases}
\]
Thus, in all cases \(|{\Fn X}| \leq 67\).
\endproof

\subsection{Triangle-free case} \label{subsec:triangle-free.case}
In this section we employ the notation and the ideas of \autoref{subsec:triangle-free.K3} and \cite[§5]{veniani1}. The following result is analogous to \cite[Proposition 5.5]{veniani1}.

\begin{proposition} \label{prop:char3-alex}
If \(X\) is triangle-free and \(D \subset \Gamma(X)\) is parabolic, then
\[
    |{\Fn X}| \leq v(D) + 24.
\]
\end{proposition}
\proof The subgraph \(D\) induces a genus 1 fibration (\cite[§3, Theorem 1]{PSS}), which can be elliptic or quasi-elliptic. The vertices in \(D \cup (\Gamma \smallsetminus \Span D)\) are fiber components of this fibration. If the fibration is elliptic, there cannot be more than 24 components, on account of the Euler--Poincaré characteristic. If the fibration is quasi-elliptic, we obtain from formula~\eqref{eq:char3-euler} that 
\begin{equation} \label{eq:char3-euler-2}
    iv^* + ii^* \leq 3.
\end{equation} 
A fiber of type~\(\IV\) can contain at most 2 lines, since there are no triangles. Hence, from \eqref{eq:char3-euler} and \eqref{eq:char3-euler-2} we deduce
\begin{align*}
    |{\Fn X}|= |\Gamma(X)| &\leq v(D) + 2\,iv + 7\,iv^* + 9\,ii^* \\
             &   = v(D) + 20 + iv^* + ii^* \\
             &\leq v(D) + 23. \qedhere
\end{align*}
\endproof

\begin{lemma} \label{lem:char3-triangle-free} 
If \(X\) is triangle-free, then \(v(\ell) \leq 12\) for any line \(\ell \subset X\).
\end{lemma}
\proof
Thanks to \autoref{lem:elliptic-triangle-free}, we can suppose that \(\ell\) is quasi-elliptic. We claim that for all \(t \in \IP\) it holds
\begin{equation} \label{eq:triangle-free2}
    2 \, v_t(\ell) \leq e(F_t) -2.
\end{equation}

Indeed, this is true if the residual cubic \(c_t\) is irreducible (\(v_t(\ell) = 0\) and \(e(F_t) \geq 2\)), or if \(c_t\) is the union of a line and a reducible conic (\(v_t(\ell) \leq 1\) and \(e(F_t) \geq 4\) because \(F_t\) is reducible).
If \(c_t\) splits into three lines (not necessarily distinct), then these lines meet in one singular point because \(X\) is triangle-free and fibers of type~\(\I_n\) do not appear in a quasi-elliptic fibration. Hence, \(F_t\) is of type~\(\IV^*\) or \(\II^*\) and \eqref{eq:triangle-free2} holds because \(v_t(\ell)\leq 3\) and \(e(F_t) \geq 8\). 

From \eqref{eq:v=sum-vF}, \eqref{eq:euler.quasi-elliptic} and \eqref{eq:triangle-free2} it follows that
\[
    2\, v(\ell) = 2\,\sum_t v_t(\ell) \leq \sum_t (e(F_t) - 2) = 24. \qedhere
\]
\endproof

\begin{proposition} \label{prop:char3-triangle-free}
If \(X\) is triangle-free, then \(|{\Fn X}| \leq 64\).
\end{proposition}
\proof
We can adapt the proof in \cite[Proposition~5.9]{veniani1}. 
Indeed, the proof of \cite[Proposition~5.8]{veniani1} is also valid in characteristic~\(3\), since we merely employ the analogous result of \autoref{prop:char3-alex} and the fact that the Picard lattice of \(X\) has rank at most \(22\). 
Therefore, if \(X\) contains no quadrangle, then \(|{\Fn X}| \leq 54\).

If \(X\) contains a quadrangle, then we can apply \autoref{lem:char3-triangle-free} to the \(\mathbf{\tilde A}_3\)-subgraph corresponding to the quadrangle, finding \(|{\Fn X}| \leq 4 \cdot (12-2) +24 = 64\).
\endproof 

\subsection{Proof of \autoref{thm:char3}} \label{subsec:proof.thm:char3}
The case in which \(X\) contains a star is treated in \autoref{prop:star}, the case in which \(X\) contains no stars, but contains a triangle is treated in \autoref{prop:triangle-not-star}, and the case in which \(X\) contains no triangles is treated in \autoref{prop:char3-triangle-free}, so the proof is now complete. \qed

\section{Examples} \label{sec:char3-examples}

In this last section, we present examples of K3 quartic surfaces with many lines defined over a field of characteristic~\(3\) (with the exception of \autoref{ex:52.lines.2.sing.pts}). Most of the examples were found during the proof of \autoref{add:char3}.
\autoref{ex:58-3rd-conf} and \autoref{ex:52.lines.2.sing.pts} were found starting from the configuration of lines communicated to the author by Degtyarev, using methods similar to the ones employed in \cite{veniani:symmetries.equations}.

\begin{example} \label{ex:58-1st-conf}
A general member of the \(1\)-dimensional family defined by
\[
x_{1}^{3} x_{2} - x_{1} x_{2}^{3} + x_{0}^{3} x_{3} - x_{0} x_{3}^{3} = a x_{0}^{2} x_{1} x_{2}
\]
is smooth and contains \(58\) lines.

More precisely, for \(a = 0\) we obtain a surface which is projectively equivalent over \(\IF_9\) to the Fermat surface and thus contains 112 lines.

If \(a \neq 0, \infty\), the surface contains a star (in \(x_0 = 0\)) formed by two cuspidal lines (of valency~30) and two elliptic lines with no other singular fibers than the star itself (hence, of valency~3). The remaining \(54\) lines are of type~\((p,q)= (1,9)\). The surface contains exactly 19 stars.
\end{example}

\begin{example} \label{ex:58-2nd-conf}
A general member of the \(1\)-dimensional family defined by
\[
x_{1}^{3} x_{2} - x_{1} x_{2}^{3} + x_{0}^{3} x_{3} - x_{0} x_{3}^{3} = a x_{0} x_{1} {\left(a x_{0} x_{2} + a x_{1} x_{3} + x_{1} x_{2} + x_{0} x_{3}\right)}
\]
is smooth and contains exactly \(58\) lines.

More precisely, as long as \(a \neq 0,\, 1,\,-1,\, \infty\), the surface contains one cuspidal line (given by \(x_0 = x_1 = 0\)) which intersects 12 lines of type~\((4,0)\), and 18 lines of type~\((1,9)\); the remaining 27 lines are of type~\((4,6)\) (for instance, \(x_2 = x_3 = 0\)). The surface contains exactly 10 stars.

For \(a = 0\) we find again a model of the Fermat surface, whereas for \(a = \pm 1\) the surface contains 20 lines and a triple point. For \(a = \infty\) we obtain the union of two planes and a quadric surface.

All surfaces of the family are endowed with the symmetries \([x_0:x_1:x_2:x_3] \mapsto [x_1:x_0:x_3:x_2]\) and \([x_0:x_1:x_2:x_3] \mapsto [x_0:x_1:-x_2:-x_3]\).
\end{example}

\begin{example}  \label{ex:58-3rd-conf}
A general member of the \(1\)-dimensional family defined by
\begin{multline*}
    {\left(a^{3} + a^{2} + a + 1\right)} \left( x_{1}^{3} x_{2} +  x_{1} x_{2}^{3} - x_{0}^{3} x_{3} - x_{0} x_{3}^{3} \right) = \\
    {\left(a - 1\right)} \left( x_{0}^{2} x_{1} x_{2} - x_{0}^{2} x_{3}^{2} + x_{1} x_{2} x_{3}^{2} \right)
    + {\left(a + 1\right)} \left( x_{1}^{2} x_{2}^{2} - x_{0} x_{1}^{2} x_{3} - x_{0} x_{2}^{2} x_{3} \right) \\
    + {\left(a^{2} - 1\right)} {\left(x_{1} x_{2} + x_{0} x_{3}\right)} {\left(x_{0} + x_{3}\right)} {\left(x_{1} + x_{2}\right)} - {\left(a^{2} + 1\right)} x_{0} x_{1} x_{2} x_{3}
\end{multline*}
is smooth and contains exactly \(58\) lines. 

More precisely, if \(a \neq 0,1,-1,\infty\) and \(a^2 \neq -1\), then the surface contains exactly one star in the plane
\begin{equation} \label{eq:char3-3rd-family-star}
x_0 + x_3 = x_1 + x_2.    
\end{equation}
The star is formed by two lines of type~\((7,0)\) and two lines of type~\((1,9)\), whose equations can be explicitly written after a change of parameter \(a = d/(d^2+1)\). Each line of type~\((7,0)\) meets 18 lines of type~\((3,6)\), and each line of type~\((1,9)\) meets 9 lines of type~\((4,6)\). All lines are elliptic.

If \(a = 0,\,1\) or \(-1\) the surface is the union of a double plane and a quadric surface. 

If \(a = \infty\) the surface is projectively equivalent to the Fermat surface.

If \(a^2 = -1\), then the surface contains \(9\) points of type~\(\bA_1\) and \(40\) lines. The star in the plane \eqref{eq:char3-3rd-family-star} is formed by two elliptic lines of type~\((4,0)\) and two quasi-elliptic lines of type~\((1,9)\). Each line of type~\((4,0)\) intersects \(9\) lines of singularity~\(2\) and valency~\(6\), while each line of type~\((1,9)\) intersects \(9\) lines of singularity~\(1\) and valency~\(9\). 

All surfaces of the family are endowed with the symmetries \([x_0:x_1:x_2:x_3] \mapsto [x_1:x_0:x_2:x_3]\) and \([x_0:x_1:x_2:x_3] \mapsto [x_0:x_1:x_3:x_2]\).
\end{example}



\begin{example} \label{ex:shimada-shioda}
The reduction modulo 3 of Shimada--Shioda's surface \(X_{56}\) (see \cite{shimada-shioda}) can be written
\[
    \Psi(x_0,x_1,x_2,x_3) = \Psi(-x_1,x_0,-x_3,x_2)
\]
where
\[
    \Psi(w,x,y,z) = w z {\left(w^{2} + w x + x^{2} + y^{2} + y z + z^{2}\right)}.
\]
It contains \(48\) lines and \(8\) singular points of type~\(\bA_1\), namely \([0:1:0:t]\), \([1:0:t:0]\), \([s:1:1:s]\) and \([s:1:-1:-s]\),  where \(t^2 + 1 = 0\) and \(s^2-s-1 = 0\). 
The plane \(\Pi_0\colon x_0 = 0\) contains a configuration \(\cB_2\) as in \autoref{fig:triangle-conf}.

The line \(x_0 = x_1 = 0\) has singularity \(0\), type~\((4,2)\) and valency \(14\). The line \(x_0 = x_2 = 0\) is quasi-elliptic and has singularity \(2\), type~\((4,6)\) and valency~\(4\).
The other two lines in \(\Pi_0\) have singularity \(1\), type~\((5,1)\) and valency~\(11\). 
There are \(8\) lines in totale through each of the two points \([0:1:0:t]\).

To the author's knowledge, this example attains the highest record of lines among non-smooth K3 quartic surfaces in characteristic~\(3\).
\end{example}

\begin{example} \label{ex:52.lines.2.sing.pts}
The complex surface defined by
\begin{multline*} 
x_{0}^2 x_{1} x_{2} - x_{1}^3 x_{2} - 2 x_{1}^2 x_{2}^2 - x_{1} x_{2}^3 + x_{0}^3 x_{3} + x_{0} x_{1}^2 x_{3} \\ + 2 x_{0} x_{1} x_{2} x_{3} + x_{0} x_{2}^2 x_{3} + 2 x_{0}^2 x_{3}^2 + x_{1} x_{2} x_{3}^2 + x_{0} x_{3}^3 = 0
\end{multline*}
contains \(52\) lines and \(2\) singular points of type~\(\bA_1\), namely \(P_0 = [0:1:-1:0]\) and \(P_1 = [1:0:0:-1]\). The surface admits the symmetries \([x_0:x_1:x_2:x_3] \mapsto [x_3:x_1:x_2:x_0]\) and \([x_0:x_1:x_2:x_3] \mapsto [x_0:x_2:x_1:x_3]\). 
The plane \(\Pi_0\colon x_0 = 0\) contains a configuration \(\cA_1\) as in \autoref{fig:triangle-conf}. Let \(\ell_1,\ell_2\) be the lines in \(\Pi_0\) of singularity \(0\) and let \(\ell_3,\ell_4\) be the lines in \(\Pi_0\) of singularity \(1\).
The lines \(\ell_1,\ell_2\) are of type~\((2,8)\) and have valency~\(14\), while \(\ell_3,\ell_4\) are of type~\((6,1)\) and have valency~\(12\). There are \(8\) lines in total through \(P_0\).

The reduction modulo \(2\) of this surface is not a K3 quartic surface because it is singular along the line defined by \(x_0 = x_1 + x_2 + x_3 = 0\).

The reduction modulo \(3\) contains \(36\) lines and \(10\) singular points of type~\(\bA_1\), namely \([0:1:-1:0]\), \([1:0:0:-1]\), \([t:1:t+1:t]\), \([t:1:1:1-t]\), \([s:1:1-s:s]\) and \([s:1:1:-1-s]\), where \(t^2-t-1 = 0\) and \(s^2 + s -1 = 0\).
The lines \(\ell_1,\ell_2\) are of type~\((2,4)\) and have valency \(10\), while \(\ell_3,\ell_4\) are of type~\((5,0)\) and have valency \(9\). There are \(6\) lines in total through \(P_0\).

The reduction modulo \(5\) contains \(56\) lines and \(4\) singular points of type~\(\bA_1\), namely \([0:1:-1:0]\), \([1:0:0:-1]\) and \([t:1:1:t]\), where \(t^2 = 3\). 
The lines \(\ell_1,\ell_2\) are of type~\((4,4)\) and have valency~\(16\), while \(\ell_3,\ell_4\) are of type~\((6,1)\) and have valency \(12\).
There are \(8\) lines in total through \(P_0\).

To the author's knowledge, this example attains the highest record of lines among non-smooth K3 quartic surfaces both in characteristic \(0\) and in characteristic \(p > 3\).
\end{example}

\bibliographystyle{amsplain}
\bibliography{references}
\end{document}